\documentclass{amsart}

\newcommand{\PSbox}[1]{\includegraphics[width=3in,height=3in]{#1}}

\usepackage{amssymb}

\usepackage[all]{xy}

\usepackage{latexsym}

\newif\ifpdfver
\ifx\undefined\pdfoutput
\usepackage[dvips]{graphicx}
\pdfverfalse
\else
\usepackage[pdftex]{graphicx}
\pdfvertrue
\fi

\input amssym.def
\input amssym.tex

\ifpdfver
\usepackage[pdftex,bookmarks=true,bookmarksopen=false,
colorlinks=true,linkcolor=blue,citecolor=blue]{hyperref}
\fi

\def\squarebox#1{\hbox to #1{\hfill\vbox to #1{\vfill}}} 
\newcommand{\stopthm}{\hfill\hfill\vbox{\hrule\hbox{\vrule\squarebox 
                 {.667em}\vrule}\hrule}\smallskip} 

\newcommand{\CC}{{\mathbb C}}

\newcommand{\RR}{{\mathbb R}}
\newcommand{\HH}{{\mathcal H}}
\newcommand{\Hh}{{\mathbb H}}

\newcommand{\SP}{{\mathbb S}}

\newcommand{\N}{{\mathbb N}}

\newcommand{\tr}{\operatorname{tr}}

\newcommand{\rest}{\!\!\restriction}

\renewcommand{\Re}{\mathop{\mathrm Re}\nolimits}
\renewcommand{\Im}{\mathop{\mathrm Im}\nolimits}

\theoremstyle{plain}

\newtheorem{prop}{Proposition}[section]
\newtheorem{lem}[prop]{Lemma}

\theoremstyle{definition}

\numberwithin{equation}{section}

\title
{
{ The Selberg zeta function for convex co-compact Schottky groups}}
\author[L. Guillop\'e]{{Laurent Guillop\'e}}
\address{Laboratoire Jean Leray  (UMR CNRS-UN 6629), D\'epartement de
  Math\'ematiques, Facult\'e des Sciences et des Techniques, 2, rue de la Houssini\`ere.
44322 Nantes Cedex 3, France}
\email{guillope@math.univ-nantes.fr}

\author[K.K. Lin]{Kevin K. Lin}
\address{Mathematics Department, University of California \\
Evans Hall, Berkeley, CA 94720, USA}
\email{kkylin@math.berkeley.edu}

\author[M. Zworski]{Maciej Zworski}
\address{Mathematics Department, University of California \\
Evans Hall, Berkeley, CA 94720, USA}
\email{zworski@math.berkeley.edu}

\begin{document}    

\begin{abstract}
  We give a new upper bound on the Selberg zeta function for a convex
  co-compact Schottky group acting on $ {\mathbb H}^{n+1}$: in strips
  parallel to the imaginary axis the zeta function is bounded by $ \exp ( C
  |s|^\delta ) $ where $ \delta $ is the dimension of the limit set of the
  group. This bound is more precise than the optimal global bound $ \exp (
  C |s|^{n+1} ) $, and it gives new bounds on the number of resonances
  (scattering poles) of $ \Gamma \backslash {\mathbb H}^{n+1} $.  The proof
  of this result is based on the application of holomorphic $
  L^2$-techniques to the study of the determinants of the Ruelle transfer
  operators and on the quasi-self-similarity of limit sets.  We also study
  this problem numerically and provide evidence that the bound may be
  optimal. Our motivation comes from molecular dynamics and we consider $
  \Gamma \backslash {\mathbb H}^{n+1} $ as the simplest model of quantum
  chaotic scattering. The proof of this result is based on the application
  of holomorphic $L^2$-techniques to the study of the determinants of the
  Ruelle transfer operators and on the quasi-self-similarity of limit
  sets.

\end{abstract}



\renewcommand{\subjclassname}{%
  \textup{2000} Mathematics Subject Classification}

\maketitle   

\section{Introduction}   
\label{in}

In this paper we give an upper bound for the Selberg zeta function of a
convex co-compact Schottky group in terms of the dimension of its limit
set. This leads to a Weyl-type upper bound for the number of zeros of
the zeta function in a strip with the number of degrees of freedom given
by the dimension of the limit set plus one. We also report on numerical
computations which indicate that our upper bound may be sharp, and close
to a possible lower bound.

Our motivation comes from the study of the distribution of quantum
resonances -- see \cite{MZ1} for a general introduction.  Since the work
of Sj\"ostrand \cite{Sj} on geometric upper bounds for the number of
resonances, it has been expected that for chaotic scattering systems the
density of resonances near the real axis can be approximately given by a
power law with the power equal to half of the dimension of the trapped
set (see \eqref{eq:conj} below). Upper bounds in geometric situations
have been obtained in \cite{WZ} and \cite{MZ}.

Recent numerical studies in the semi-classical and several convex
obstacles settings, \cite{KL},\cite{LZ} and \cite{lsz} respectively,
have provided evidence that the density of resonances satisfies a lower
bound related to the dimension of the trapped set. In complicated
situations which were studied numerically, the dimension is a delicate
concept and it may be that different notions of dimension have to be
used for upper and lower bounds -- this point has been emphasized in
\cite{lsz}.

Generally, the zeros of dynamical zeta functions
are  interpreted as the classical correlation
spectrum  \cite{R}.
In the case of {\em convex co-compact hyperbolic quotients}, $ X = \Gamma 
\backslash \Hh^{n+1} $ 
quantum resonances also 
coincide with the zeros of the zeta function -- see
\cite{PP}. 
The notion of the dimension of the trapped set is also clear as it
is given by $ 2 ( 1 + \delta ) $. Here $ \delta = \dim \Lambda( \Gamma ) $ 
is the dimension of the limit set of $ \Gamma$, that is the set of 
accumulation points of any $\Gamma$-orbit in $\Hh^{n+1}$,
$ \Lambda ( \Gamma ) \subset \partial \Hh^{n+1} $. 

Hence we may expect that
\begin{equation}
\label{eq:conj}
  \sum_{ |\Im s | \leq r \,, \; \Re s > - C } m_\Gamma
 ( s ) \sim r^{ 1+ \delta}
\, , \end{equation} 
where $ m_\Gamma
 ( s ) $ is the multiplicity of the zero of the zeta function of
$ \Gamma $ at $ s $. 

Referring for definitions of Schottky groups and zeta functions to 
Sections \ref{sg} and \ref{pz} respectively we have 

\medskip
\noindent
{\bf Theorem.} {\em Suppose that $ \Gamma $ is a convex co-compact 
Schottky group and that $ Z_\Gamma  ( s ) $ is its Selberg zeta function.
Then for any $ C_0 > 0 $ there exists $ C_1 $ such that for 
$ | \Re s |< C_0 $
\begin{equation}
\label{eq:zetad}
  | Z_\Gamma ( s ) | \leq C_1 \exp ( C_1 |s|^{\delta } ) \,, \ \ 
\delta = {\mathrm{dim}} \; \Lambda ( \Gamma ) \,. \end{equation}}

\medskip

The proof of this result is based on the {\em quasi-self-similarity} 
of limit sets of convex co-compact Schottky groups and on  
the application of holomorphic $ L^2$-techniques
to the study of the determinants of the Ruelle transfer operators.

\begin{figure}[htb]
\PSbox{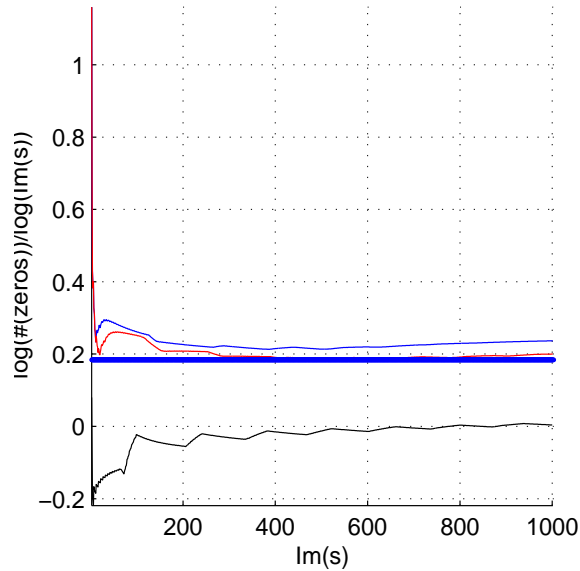}
\caption{The plot of $ \log\log N(k) / \log k - 1 $, where 
$ N(k ) $ is 
the number of of zeros with $ |\Im s| \leq k $,
for a Schottky reflection group with $ \delta \simeq
0.184 $. Different lines represent different strips, and the thick 
blue line gives $ \delta$.}
\label{f30}
\end{figure}

If we use the convergence of the product representation
\eqref{eq:zeta}  of the zeta function for $ \Re s $ large and apply 
Jensen's theorem we obtain the following 

\medskip
\noindent
{\bf Corollary 1.} {\em Let $ m_\Gamma (s) $ be the multiplicity of a zero
of $ Z_\Gamma $ at $ s$. Then, for any $C_0$, there exists some constant
$C_1$ such that for $ r > 1 $
\begin{equation}
\label{eq:ubd}
  \sum \{ 
m_\Gamma (s) \; : \;  {r \leq |\Im s | \leq r + 1}\,, \ \ 
 { \Re s > - C_0 } 
\}  \leq  C_1 r^{ \delta}
\, , \end{equation} 
where $ \delta = {\mathrm{dim}} \; \Lambda ( \Gamma ) $. }

\medskip

We can apply the preceding results to Schottky manifolds: a hyperbolic
manifold is said Schottky if its fundamental group is Schottky. The case of
surfaces being of special interest: any convex co-compact hyperbolic
surface is Schottky. With the description of the divisor of the zeta
function through spectral data established by Patterson and Perry
\cite{PP}, we can reformulate the preceding corollary nicely in the
resonance setting. We do it only for surfaces (see below for short comments
on higher dimensions)

\medskip
\noindent
{\bf Corollary 2.} {\em Let $M$ be a convex co-compact hyperbolic surface,
  $\mathcal{S}_M$ be set of the scattering resonances of the
  Laplace-Beltrami operator on $M$ and $ m_M (s) $ be the multiplicity of
  resonance $s$. Then, for any $C_0$, there exists some constant $C_1$
  such that for $ r > 1 $
\begin{equation}
\label{eq:Surf}
  \sum \{ 
m_M (s) \; : \; s\in\mathcal{R}, {r \leq |\Im s | \leq r + 1}\,, \ \ 
 { \Re s > - C_0 } 
\}  \leq  C_1 r^{ \delta}
\, , \end{equation} 
where $ 2(1+\delta)$ is the Hausdorff dimension of the 
recurrent set for the geodesic flow on $T^*M$.}
\medskip 

This corollary is stronger than the result obtained 
in \cite{MZ}
where the upper bound of the type \eqref{eq:conj} was given. In fact, the
upper bound \eqref{eq:Surf} is what we would obtain had we had a Weyl law of
the form $ r^{ 1 + \delta } $ with a remainder $ {\mathcal O} ( r^{\delta }
) $. That local upper bounds of this type are expected despite the absence
of a Weyl law has been known since \cite{PZ}.

Section \ref{num} deals with numerical computations of the density of 
zeros. They show that \eqref{eq:conj} may be true\footnote{Strictly
speaking, the numerical computations are done for a slighly different
zeta function of conformal dynamical system}. In fact, in 
the range of $ \Im s $ used in the computation 
we see that the number of zeros grows fast. If the range of  $ \Re s $ 
is large (and fixed) we need very large $ \Im s $ to see the upper bound
of Corollary 2. The computations also show that our bound on the zeta
function is optimal. For values of $ Z_\Gamma ( s ) $ with $ \Re s $ negative
we see that we need very large $ \Im s $ to see the onset of the upper
bound. That is not surprising since we recall in Proposition \ref{p:1}
 that $ \log | Z_\Gamma ( s) | 
= {\mathcal O} ( |s|^{n+1} ) $, and that this bound is optimal (and
of course $ \delta < n $). 

\begin{figure}[htb]
\begin{center}
\PSbox{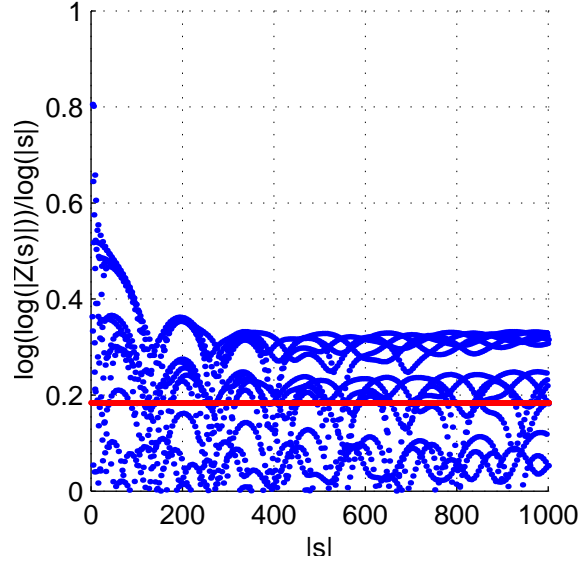}
\end{center}
\caption{Density of values of $ \log|\log|Z_{\Gamma } ||/\log|s|  $
for a Schottky reflection group with $ \delta \simeq
0.184 $.}
\label{fd30}
\end{figure}

\begin{figure}[htb]
\includegraphics[width=3in]{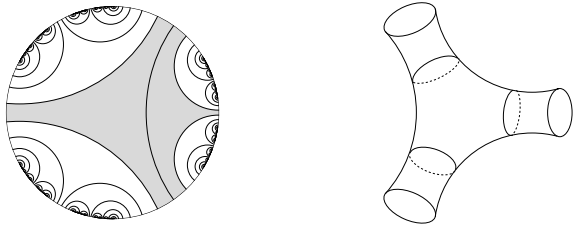}
\caption{Tesselation by the group, $ \Gamma_\theta $, 
$ \theta = 110^\circ $, 
generated by symmetries in three symmetrically placed lines each cutting the
unit circle in an 110$^\circ$ angle, with the fundamental domain of its
Schottky subgroup of direct isometries, $ \Gamma^+_\theta $, and the
associated Riemann surface $ \Gamma^+_\theta \backslash \Hh^2 $. The
dimension of the limit set is $\delta=0.70055063\ldots$.}
\label{fi:0}
\end{figure}

We refer to Section \ref{num} for the details and present here two
pictures only. 
We take 
for $ \Gamma $ a group generated by compositions of reflections in three
symmetrically spaced circles perpendicular to the unit circle, and cutting
it at the angles $ 30^\circ $ (see Fig.\ref{fi:0} for the $ 110 ^\circ $ angle). 
Fig.\ref{f30} shows the density of zeros of $ Z_\Gamma $ in that case and
Fig.\ref{fd30} plots the values of $ \log|\log|Z_{\Gamma } || $. 

\section{Schottky groups}
\label{sg}

The hyperbolic geometry on $\Hh^{n+1}$ and the conformal geometry on its
boundary at infinity $\partial \Hh^{n+1} =\SP^n$ share the same automorphism group:
the isometry group $\mathrm{Isom} ( \Hh^{n+1}) $ and the conformal group
$\mathrm{Conf}(\SP^n)$ (with the conformal structure given by the standard
metric on $ \SP^n $ of curvature $+1$) are isomorphic.  In particular any
isometry $g$ of $ \Hh^{n+1}$ induces on $ \SP^n $
a conformal map $\gamma$, whose conformal distorsion at the point
$w\in\SP^n$ will be denoted by $\|D\gamma(w)\|$. There is also a
correspondence between balls $\mathcal{D}$ and spheres $\mathcal{C}$ on
$\SP^n$ (for $n=2$, the original setting for Kleinian groups, these are
discs and circles) and half-spaces $\mathcal{P}$ and geodesic hyperplanes
$\mathcal{H}$ in $\Hh^{n+1}$: $\mathcal{D}=\overline{\mathcal{P}}\cap
\partial\Hh^{n+1}$ and $\mathcal{C}=\overline{\mathcal{H}}\cap
\partial\Hh^{n+1}$. Given a hyperplane $\mathcal{H}$ (a sphere
$\mathcal{C}$ resp.), its interior hyperplane (ball) will be given by the
choice of a component of $\Hh^{n+1}\setminus\mathcal{H}$
($\SP^n\setminus\mathcal{C}$ resp.)

Let us review  definitions of a Schottky group (see \cite{McM1}, \cite{Mas2}, 
\cite{Rat} and references given there), defined originally in
conformal geometry on $S^2$.
To the configuration of mutually disjoint geometric balls  on the sphere $S^n$,
$\mathcal{D}_i,i=1,\ldots\ell$, 
we associate the
Schottky marked reflection group,
\[ \Gamma(\mathcal{D}_1,\ldots,\mathcal{D}_\ell)\,,\]
defined as the group generated by
the inversions $\sigma_i$ in $\partial\mathcal{D}_i,i=1,\ldots,\ell$.  The
corresponding hyperbolic group is the Schottky marked reflection group
\[ \Gamma(\mathcal{P}_1,\ldots,\mathcal{P}_\ell)\,,\] 
for the  collection,
$\mathcal{P}_1,\ldots,\mathcal{P}_\ell$, of mutually disjoint hyperbolic
half-spaces $\mathcal{P}_i$ with boundaries at infinity given by 
 $\mathcal{D}_i$'s.  The
group is generated by the hyperbolic symmetries $s_i$ in the hyperplane
$\partial\mathcal{P}_i$ (with infinite boundary $\partial\mathcal{D}_i$),
$i=1,\ldots\ell$.

More generally, to a configuration of $2\ell$ disjoint topological balls
$$\mathcal{D}_i\,, \ \ i=1,\ldots,2\ell\,, $$
with a set of conformal maps 
$\gamma_i,i=1,\ldots,\ell$, satisfying
\[ \gamma_i(\SP^n\setminus\mathcal{D}_i)=\overline{\mathcal{D}_{\ell+i}}\,,\]
we 
associated the Schottky marked group
\[ \Gamma(\mathcal{D}_1,\ldots,\mathcal{D}_{2\ell},\gamma_1,\ldots,
\gamma_\ell)\,,\]
generated by the $\gamma_i,i=1,\ldots,\ell$. 

We also introduce maps
\[ \gamma_{\ell+j}=\gamma_j^{-1}, j=1,\ldots,\ell \,, \]
which map the exteriors
of $\mathcal{D}_{\ell+j}$'s onto $\overline {\mathcal{D}_j}$. 

The Schottky group $ \Gamma $ is a free group on the $\ell$ (free) generators
$\gamma_1,\dots,\gamma_\ell$. If the topological balls
$\mathcal{D}_i,i=1,\ldots,2\ell$ are geometric balls, the marked Schottky
group is said to be \emph{classical}. A classical Schottky marked group
$\Gamma(\mathcal{D}_1,\ldots,\mathcal{D}_{2\ell},\gamma_1,\ldots,\gamma_\ell)$,
subgroup $\mathrm{Conf}(\SP^n)$, can be presented as an isometry group of the
hyperbolic space $\Hh^{n+1}$. 

Let $\mathcal{P}_i$ be the hyperbolic
half-space with boundary at infinity given by $\mathcal{D}_i$.
If $g_i$ is the
hyperbolic isometry of $\Hh^{n+1}$ with action at infinity given by 
$\gamma_i$, then the
classical Schottky group
\[ \Gamma=\Gamma(\mathcal{D}_1,\ldots,\mathcal{D}_{2\ell},\gamma_1,\ldots,\gamma_\ell)\,,\]
has a hyperbolic marking 
\[ \Gamma=\Gamma(\mathcal{P}_1,\ldots,\mathcal{P}_{2\ell},g_1,\ldots,
g_\ell)\,.\]
The Schottky domain $\Hh^{n+1}\setminus\cup_{i=1}^{2\ell}\mathcal{P}_i$
is a fundamental domain for the action of $\Gamma$ on $\Hh^{n+1}$.

A group is said to be a Schottky (reflection) group if it admits a presentation
induced by a 
configuration of balls as described above above. We remark that
the subgroup 
of positive isometries of a Schottky reflection
group is a classical Schottky group, with a marking
\[ (\mathcal{D}_1,\ldots\mathcal{D}_{\ell-1},
\sigma_\ell\mathcal{D}_1,\ldots, \sigma_\ell\mathcal{D}_{\ell-1},
g_1,\ldots, g_{\ell-1})\,, \]
where $g_i=\sigma_{\ell}\sigma_i $.   Such
a group was called symmetrical by Poincar\'e \cite{Poin}. 

 An oriented
hyperbolic manifold $M$ is said to be (classical) Schottky if its fundamental
group $\pi_1(M)$ (realized as a discrete subgroup of
$\mathrm{Isom}^+(\Hh^{n+1})$) admits a (classical) Schottky marking.

A Schottky group is convex co-compact if in addition
the closures $\overline{\mathcal{D}_i}$ are mutually disjoint. In the 
convex co-compact case all 
non elliptic elements of $ \Gamma $ are 
hyperbolic. That means in particular that for $ \gamma \in \Gamma $
there exists $ \alpha \in  {\mathrm{Isom}} ( \Hh^{n+1}) $ such that, in the
Poincar\'e model $\Hh^{n+1} \simeq \RR^{n+1}_+= \RR_+\times \RR^n $,
\begin{equation}
\label{eq:hyp}
  \alpha^{-1} \gamma \alpha ( x ,y ) =  e^{ \ell ( \gamma ) } (x, 
O(\gamma) y )  \,,
\ \ (x, y ) \in \RR^{n+1}_+\,,
\ \  O(\gamma) \in O(n ) \,, \ \ \ell ( \gamma ) > 0 \,. 
\end{equation}
If $\Gamma\subset\mathrm{Isom}^+(\Hh^{n+1})$, the conjugacy classes of hyperbolic elements,
\[ \{ \gamma_1 \} = \{ \gamma_2 \}  \ \Longleftrightarrow \ 
\exists \; \beta \in \Gamma \  \beta \gamma_1 \beta^{-1} = \gamma_2 \,,\]
are in one-to-one correspondence with 
closed geodesics of $ \Gamma \backslash \Hh^{n+1} $. The primitive 
geodesics correspond to conjugacy classes of primitive elements of 
$ \Gamma $ (that is, elements which are not non-trivial powers).
The magnification factor $\exp( \ell ( \gamma )) $ in \eqref{eq:hyp}
gives the length $\ell(\gamma)$ of the closed geodesic. 

The limit set, $ \Lambda ( \Gamma ) $  of a discrete subgroup, $ \Gamma$,
of $ {\mathrm{Isom}}( \Hh^{n+1} ) $, is defined as the set 
in $\overline{\Hh^{n+1}}=\Hh^{n+1}\cup \partial \Hh^{n+1}$ of accumulation
points of any $\Gamma$-orbit in $\Hh^{n+1}$: the limit set $ \Lambda (
\Gamma ) $ is included in the boundary $ \partial\Hh^{n+1}$. In the convex
co-compact case it has a particularly nice structure; furthermore, for
Schottky groups, it is totally disconnected and included in
$\mathcal{D}=\cup_{i=1}^{2\ell}\mathcal{D}_i$. The aspects relevant to us come
  from the work of Patterson and Sullivan -- see \cite{Sull} and references
  given there. As will be discussed in more detail in Section \ref{mp}, the
  limit set has a quasi-self-similar structure and a finite Hausdorff
  measure at dimension $ \delta = \delta ( \Gamma ) $.

The limit set is related to the {\em trapped set}, $ K $,  of the usual
scattering \cite{Sj}, \cite{WZ}, that is the set of points in 
phase space such that the trajectory through that point does not 
escape to infinity in either direction.: if $\Delta$ is the diagonal of
$K\times K$ and $\pi^*$ the projection from $T^*\Hh^{n+1}$ on
$T^*\Gamma\backslash\Hh^{n+1}$,  the trapped set $K$ is the union of the projections
$\pi^*(C_{\xi\eta})$ where $C_{\xi\eta}$ is the geodesic with extremities
$\xi$ and $\eta$, both in the limit set $\Lambda(\Gamma)$. In particular,
 we have 
\[ \dim K = 2 (\delta + 1)\,, \]
see \cite{MZ}.

To stress the connection to closed orbits let us also mention that,
generalizing earlier
results of Guillop\'e \cite{G} and Lalley \cite{L}, Perry \cite{Pe2} showed that
\[ \sharp \{ \{ \gamma \} \;: \; \gamma \ \text{primitive,} \
\ell ( \gamma ) < r \} \sim 
\frac{ e^ {\delta r }} {\delta r } \,. \]

\begin{figure}
\begin{center}
\includegraphics[width=3in]{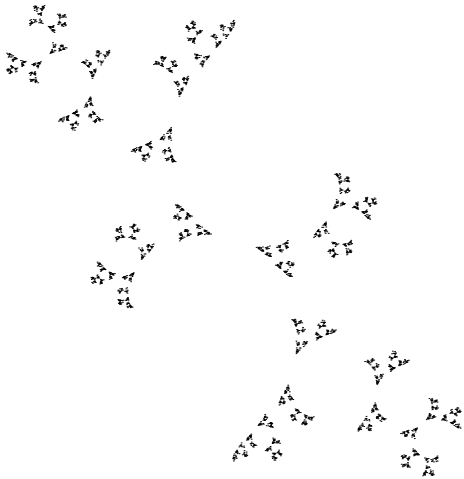}
\end{center}
\caption{A typical limit set for a convex co-compact 
Schottky group, $ \Gamma  
\subset {\mathrm{Isom}} (\Hh^3) $.}
\label{fi:2}
\end{figure}

\section{Properties of the Selberg zeta function}
\label{pz}

For $ \Gamma $, a discrete subgroup of $ {\mathrm{Isom}} (\Hh^{n+1}) $,
the Selberg zeta function is defined as follows
\begin{equation}
\label{eq:zeta}
 Z_\Gamma ( s ) =  \prod_{ \{ \gamma \} } \prod_{\alpha \in \N_0^n} 
\left( 1 -  e^{ - i \langle \theta( \gamma ) , \alpha \rangle }
e^{ - ( s + |\alpha| )  \ell( \gamma ) } \right) \,.
\end{equation}
Here, $ \gamma \in \Gamma  $ are hyperbolic, $ \exp ( \ell ( \gamma ) 
+ i \theta_j ( \gamma ) ) $ are the eigenvalues of the derivative of
the action of $ \gamma $ on $ \SP^n $ at the repelling fixed point
$(\exp(i\theta_j(\gamma))$ are the eigenvalues of the isometry $O(\gamma)$
in the normal form \eqref{eq:hyp}) and
\[ \{ \gamma \} = \text{ the conjugacy class of a primitive hyperbolic 
element $ \gamma $.}\]
An element is called primitive if it is {\em not} a non-trivial power
of another element. 

We define 
the following map on $ \mathcal{D} = \bigcup_{ i=1}^{2\ell}  \mathcal{D}_i$:
\begin{equation}
\label{eq:T}
 T: \mathcal{D} \; \longrightarrow \; \SP^n \,, \ \ T (x) = g_i ( x )  \,, 
\ \ x \in \mathcal{D}_i \,.\end{equation}
We need to find an open neigbourhood of the limit set where $T$ is
strictly expanding in the following sense~: $F$ defined on $V$ is
said to be 
(strictly) expanding on $V$ with respect to the metric $\|\hphantom{\xi}\|$
if there exists $\theta\ge 1$ ($\theta>1$) such that
\[
\|DF(v)\xi\|\ge \theta\|\xi\|\,,\quad v\in V,\xi\in T_vV\,.
\]

In the case when $\Gamma$ is the positive part of a Schottky reflection
group, we can suppose that up to a conformal identification
$\partial\mathcal{D}_\ell$ is a great circle of the sphere $S^n$.
For  the metric on $S^n$ we can take the metric
induced by its embedding in $\RR^{n+1}$. The inversion,
$\sigma_\ell$, is the restriction to $S^n$ of the symmetry on $\RR^{n+1}$
with respect of the euclidean hyperplane countaining 
$\partial\mathcal{D}_\ell$,
hence an isometry. Each inversion $\sigma_i,i=1,\ldots,\ell-1$ is
expanding on the ball $\mathcal{D}_i$. Hence, the map $T$ is expanding on
$V=\bigcup_{i=1}^{\ell-1}\mathcal{D}_i\cup\sigma_\ell\mathcal{D}_i$,
strictly expanding on any open set precompact in $V$.

However, the map $T$ is not expanding on $\mathcal{D}$ in general.
To circumvent that we need to consider refinements,
$\mathcal{D}^N$, defined by recurrence:
\[ \begin{split}\mathcal{D}^1 & =\mathcal{D}\\
\mathcal{D}^N & 
=T^{-1}(\mathcal{D}^{N-1})\cap \mathcal{D}^{N-1}\,.
\end{split} \] 
Each set
$\mathcal{D}^N$ is a disjoint  union $\cup_{i=1}^{d_N}\mathcal{D}^N_i$.
The collection of sets $\{ \mathcal{D}^N_i\}_{i} $ coincides with the 
collection
\[  \{ \mathcal{D}_\gamma\}_{\gamma \in 
\Gamma_N }\,, \ \ \Gamma_N = \{ \gamma \in \Gamma \; : \; |\gamma| =N \} \,.\]
Here  $|\gamma|$ is the 
combinatorial length with respect to the system of generators
$\{g_1,\ldots,g_{2\ell}\}$. In other words,
\[
\mathcal{D}_\gamma=x_N^{-1}\ldots x_{2}^{-1}\mathcal{D}_{x_1}
 \quad \hbox{ if }\quad
\gamma=x_1\ldots x_N\,,
\]
where $\mathcal{D}_{x_N}=\mathcal{D}_i$ if $x_N=g_i$.  The iterated map, 
$T^{N}$,  is defined on $\mathcal{D}^{N}$, and
\[ T^{|\gamma|}\mathstrut_{|\mathcal{D}_\gamma}=\gamma\,, \]
The map $T$ is
strictly expanding on $\mathcal{D}^N$ for $N$ big enough as explained in the
following lemma (see Lemma 9.2 in Lalley
\cite{L} for a similar result).
\begin{lem}
  Let $\Gamma$ a Schottky group and $\mathcal{D},T$ defined as in
  \eqref{eq:T}. There exist an integer $N\ge 1$, a metric $\|\hskip
  1em\|_\Gamma$, defined on $\mathcal{D}^N$, and a real $\beta>1$ such that
  $\|DT(w)\|_\Gamma\ge \beta, w\in \mathcal{D}^N$. The metric can be taken
  analytic on $\mathcal{D}^N$.
\end{lem}
\begin{proof}
  
  Let us recall that any M\"obius transformation $\gamma$ of $\RR^{k}$
  which does not
  fix the point at infinity, $\infty$, has an isometric sphere $S_\gamma$, and
that  $\gamma^{-1}$ is strictly contracting  on any compact subset of its exterior
  (the unbounded component of
  $\RR^k\setminus S_\gamma$). The sphere $S_\gamma$
  is centered at $\gamma^{-1}\infty$ and if $r_\gamma$  is its radius, we
  have (see \cite{Rat})
  \begin{equation}
    \label{eq:MobDis}
    ||\gamma x-\gamma
      y||=\frac{r_\gamma^2\,||x-y||}{||x-\gamma^{-1}\infty||\,
      ||y-\gamma^{-1}\infty||}\,,\quad
x,y\in \RR^k\setminus\{\gamma^{-1}\infty\}
  \end{equation}
  
Up to a conformal transformation, we can suppose that $\Gamma$ is a subset of
$\mathrm{Conf}(\RR^n)$, with the point at infinity in its ordinary set.
No non trivial element in $\Gamma$ fixes $\infty$ and, taking in
\eqref{eq:MobDis} as $x,y$ points in the upper half-plane $\Hh^{n+1}$ (and
the Poincar\'e extension of $\gamma$ to $\RR^{n+1}$), we deduce that the 
set of radii $\{r_\gamma,\gamma\in\Gamma\}$ accumulates only at $0$.

For $\gamma=x_1\ldots x_N$, we have
\[
\gamma^{-1}(\infty)\in{\mathcal{D}}_{\gamma}\subset{\mathcal{D}}_{
  x_2\ldots  x_N}\subset\ldots\subset
\overline{{\mathcal{D}}_{ x_{N-1}
  x_N}}\subset{\mathcal{D}}_{ x_N}
\]
hence there exists $N_0$ such that the isometric sphere
$S_{\gamma}$ is included in ${\mathcal{D}}_{ x_N}$ if
$|\gamma|\ge N_0$.  For such a $\gamma$, the interior of the
isometric sphere $S_{\gamma^{-1}}$ is included in
${\mathcal{D}}_{ x_1^{-1}}$: its exterior contains all the
$\overline{\mathcal{D}_{ x_i}}, x_i\not=x_1^{-1}$, hence
$\gamma$ is strictly contracting on 
$\cup_{x_i\not=x_1^{-1}}\mathcal{D}_{ x_i}$.  As $
\gamma({\mathcal{D}}_{\gamma})\subset\cup_{x_i\not=x_1^{-1}}{\mathcal{D}}_{
  x_i}$, the map $\gamma=
T^{|\gamma|}\mathstrut_{|\mathcal{D}_{\overline\gamma}}$ is expanding on
$\mathcal{D}_\gamma\subset{\mathcal{D}}^{|\gamma|}$.  We have
just proved the existence of $\eta_0>1$ such that
$\|DT^{N_0}(w)\|\ge\eta_0, w\in\mathcal{D}^{N_0}$, hence there exist
constants $C>0, \theta>1$ such that $\|DT^p(w)\|\ge C\theta^p, w\in
\mathcal{D}^p, p\ge 1$.

Taking an integer $N$ such that $C\theta^N>1$, we define on $\mathcal{D}^N$
the metric (introduced by Mather \cite{M})
\[
\|V\|_\Gamma=\sum_{p=0}^{N-1}\|DT^p(w)V\|\,,\quad V\in T_w\mathcal{D}^N\,,
\]
which concludes the proof.
\end{proof}

Let $ \widetilde \SP^n $ be a  Grauert tube of $ \SP^n $, that is
a complex $n$-manifold containing $ \SP^n $ as a totally real 
submanifold (that is all we need). Let us then choose 
open neighbourhoods\footnote{we drop the index $N$ in the
open sets $D_i^N $ for the purpose of  notational simplicity} 
 $D_i, i=1,\ldots,d_N$ of $\mathcal{D}_i^N$ in
$\widetilde\SP^n$.  By further shrinking, we can suppose that the open sets
$D_i$ are mutually disjoint, and that the real analytic maps $T$ and
$\|DT\|_\Gamma$ extend holomorphically to $D=\cup_{i=1}^{d_N}D_{i}$, with
$\|DT\|_\Gamma\ge\widetilde\beta$ for some $\widetilde\beta>1$. The open
sets $D_i$ can be chosen to be
a union $D_i=\cup_{k=1}^{\delta_i}D_{ik}$ of open
sets, each one biholomorphic to the ball $B_{\CC^n}(0,1)$ in $\CC^n$.

With this formalism in place we define the Ruelle transfer operator
\begin{gather}
\label{eq:trans}
\begin{gathered}
{\mathcal L} ( s ) u ( z ) = \sum_{ Tw = z } \| DT ( w ) \|^{-s} 
u ( w ) \,,  \ \ z \in D \,, \\
u \in H^2 ( D) \,, \ \ H^2 ( D) = \{ u \text{ holomorphic in $ D$ } \;: \; 
 \int \! \!\! \int_D |u( z ) |^2  dm ( z) < \infty \} \,.
\end{gathered}
\end{gather}
The only difference from the standard definition 
lies in choosing $ L^2 $ spaces of holomorphic 
functions instead of Banach spaces. However we still obtain the 
analogue of a (special case of a) 
result of Ruelle \cite{Rue} and Fried \cite{fr}:

\begin{prop}
\label{p:1}
Suppose that $ {\mathcal L} ( s): H^2 ( D) \rightarrow H^2 ( D) $ is
defined by \eqref{eq:trans}. Then for all $ s \in \CC $ $ {\mathcal L } ( s) $
is a trace class operator and 
\begin{equation}
\label{eq:det}
| \det ( I - {\mathcal L} ( s ) ) | \leq \exp ( C |s|^{n+1} ) \,.
\end{equation}
\end{prop}
\begin{proof}
The proof is based on estimates of the characteristic values, $ \mu_\ell (
{\mathcal L} ( s ) ) $. We will show that there exists $ C > 0 $ such that
\begin{equation}
\label{eq:cest}
 \mu_\ell ( {\mathcal L}( s ) ) \leq Ce^{ C|s| - \ell^{\frac{1}{n}}/ C } \,.
\end{equation}

To see how that is obtained and how it implies \eqref{eq:det} let us
first recall some basic properties of  characteristic values of
a compact operator $ A: H_1 \rightarrow H_2 $ where $ H_j $'s are Hilbert
spaces. We define
\[ \| A \| = \mu_0 ( A ) \geq \mu_1 ( A ) \geq \cdots \geq \mu_{\ell} ( A ) 
\rightarrow 0 \,, \]
to be the eigenvalues of $ ( A^* A )^{\frac12}: H_1 
\rightarrow H_1 $, or equivalently of $ ( A A^*)^{\frac12}: H_2 \rightarrow
H_2 $. The min-max principle shows that
\begin{equation}
\label{eq:mm} 
\mu_\ell ( A ) = \min_{ {V \subset H_1} \atop
{{\mathrm{codim}}\; V = \ell}} \max_{ {v\in V}\atop{ \|v\|_{H_1} = 1}} 
\| A v\|_{H_2}  \,.
\end{equation}
The following rough estimate will be enough for us here: suppose that 
$ \{ x_j\}_{ j =0}^\infty $ is an orthonormal basis of $ H_1 $, then 
\begin{equation}
\label{eq:cec} 
\mu_{\ell} ( A ) \leq \sum_{ j = \ell}^\infty \| A x_j \|_{H_2} \,.
\end{equation}
To see this we will use $ V_\ell = {\mathrm{span}}\; \{ x_j \}_{j=\ell}^\infty$
in \eqref{eq:mm}: for $ v \in V_\ell $ we have, by the Cauchy-Schwartz
inequality,  and the obvious $ \ell^2 \subset \ell^1 $ inequality,
\[ \| A v \|_{H_2}^2 = 
\left\| \sum_{ j=\ell}^\infty \langle v, x_j \rangle_{H_1} A x_j \right\|
\leq 
\| v \|_{ H_1}^2 \left( \sum_{ j = \ell}^\infty \| A x_j \|_{H_2} \right)^2 
\,,\]
from which \eqref{eq:mm} gives \eqref{eq:cec}.

We will also need some real results about characteristic values
The first is the {\em Weyl inequality} (see \cite{GK}, and also
\cite[Appendix A]{Sj}). It says that if $ H_1 = H_2 $ and
$ \lambda_j ( A ) $ are the eigenvalues of $ A$, 
$$ |\lambda_0 ( A) | \geq
|\lambda_1 ( A) | \geq \cdots \geq |\lambda_\ell ( A ) |\rightarrow 0 \, ,$$ 
then for any $ N $, 
\[ \prod_{ \ell=0}^N ( 1 + |\lambda_\ell ( A ) | ) \leq 
\prod_{ \ell=0}^N ( 1 + |\mu _\ell ( A ) | ) \,.\]
In particular if the operator $ A $ is of trace class, that is if,
$ \sum_{\ell} \mu_\ell ( A ) < \infty $, then the determinant
\[ \det ( I + A )  \stackrel{\mathrm{def}}{=} \prod_{ \ell=0}^\infty ( 1 + 
\lambda_\ell ( A ) ) \,,\]
is well defined and
\begin{equation}
\label{eq:weyl}
 |\det( I + A ) | \leq 
\prod_{ \ell=0}^\infty ( 1 + \mu_{\ell} ( A ) ) \,.
\end{equation}

We also need to recall the following standard inequality about 
characteristic values (see \cite{GK}):
\begin{equation}
\label{eq:stand}
 \mu_{\ell_1 + \ell_2 } ( A + B ) = \mu_{\ell_1} ( A ) + \mu_{\ell_2} ( B) 
\end{equation}
We finish the review, as we started, with an obvious equality:
suppose that $ A_j: H_{1j} \rightarrow H_{2j} $ and we form 
$ \bigoplus_{j=1}^J A_j:  \bigoplus_{j=1}^J H_{1j} \rightarrow  
\bigoplus_{j=1}^J 
H_{2j} $, as usual, $ \bigoplus_{j=1}^J A_j ( v_1 \oplus \cdots \oplus v_J ) =
A_1 v_1 \oplus \cdots \oplus A_J v_J $. Then
\begin{equation}
\label{eq:obv}
\sum_{\ell=0}^\infty 
\mu_\ell 
 \left(\bigoplus_{j=1}^J A_j \right) = \sum_{j=1}^J \sum_{\ell=0}^\infty 
\mu_\ell ( A_j ) \,. 
\end{equation}

With these preliminary facts taken care of, we see that \eqref{eq:cest}
implies \eqref{eq:det}. In fact, \eqref{eq:weyl} shows that
\[ |\det ( I - {\mathcal L} ( s) )| \leq \prod_{ \ell=0}^\infty  (  1 + e^{C|s| 
- \ell^{\frac{1}n} / C} ) \leq e^{ C_1 |s|^{n+1} } \,.\]
Hence it remains to establish \eqref{eq:cest}. For that we will 
write 
\[ H^2 ( D) = \bigoplus_{i=1}^{d_N} H^2 ( D_i ) \,,\]
and introduce, for $i,j=1,\ldots,d^N$, the operator 
\[
 {\mathcal L}_{ i j } ( s ): H^2 ( D_i ) \rightarrow H^2 ( D_j ) \,,
\]
non zero only when $T(D_i)$ and $D_j$ are not disjoint, where
\begin{equation}
  \label{eq:fij}
{\mathcal L}_{ij} ( s ) u ( z) \stackrel{\mathrm{def}}{=}  \| Df_{ij} ( z)\|^s 
u ( f_{ i j } ( z ) ) \,, \ \ z \in D_j \,, \ 
\ f_{ij} = (T \rest_{ D_i} )^{-1} \rest_{ D_ j} \,. 
\end{equation}
From \eqref{eq:stand} and a version of \eqref{eq:obv} we then have
\[ \mu_\ell ( {\mathcal L} ( s) ) \leq \max_{ {1 \leq i, j \leq d_N }
  } 2 \mu_{ [\ell/ C] } ( {\mathcal L}_{ij} ( s) ) \,. \] 

To estimate $ \mu_k ( {\mathcal L}_{ij} (s)) $, let us recall that 
$D_i$ was taken as an
union of open sets $D_{ik}, k=1,\ldots,\delta_i$ biholomorphic to $B_{\CC^n} (
0, 1 )$: as $f_{ij}(D_j)$ is relatively compact in $D_i$, we can find
$\rho\in(0,1)$ (independent of $i,j=1,\ldots,d_N$) such that
$f_{ij}(D_j)\subset D_i^\rho$ where
$D_i^{\rho}=\cup_{k=1}^{\delta_i}D_{ik}^\rho$ with $D_{ik}^{\rho}\subset
D_{ik}$ the pullback of the ball $B_{\CC^n} ( 0, \rho )$ through the
biholomorphism of $D_{ik}$ onto $B_{\CC^n} ( 0, 1 )$. The map
$\mathcal{L}_{ij} (s)$ is the composition
\[
\xymatrix{H^2(D_i)\ar[r]^-{R}&\oplus_{k=1}^{\delta_i}
  H^2(D_{ik})\ar[r]^-{\oplus R^\rho_{ik}}&
\oplus_{k=1}^{\delta_i} H^2(D_{ik}^\rho)\ar[r]^-{\pi}&
H^2(D_i^{\rho})\ar[r]^-{\mathcal{L}^{\rho}_{ij}(s)}& H^2(D_j)}
\]
where $R$ and $R^\rho_{ik}$ are the natural restrictions, $\pi^\rho$ is the
orthogonal projection on the space $H^2(D_i^{\rho})$ immersed in $\oplus_{k=1}^{\delta_i}
H^2(D_{ik}^\rho)$ by the natural restrictions and
$\mathcal{L}^{\rho}_{ij}(s)$ is defined by the same formula
\eqref{eq:fij} as
$\mathcal{L}_{ij}(s)$.  The maps $R$ and $\pi^\rho$ are bounded,
while the norm of $\mathcal{L}_{ij}(s)$ is bounded by $Ce^{C|s|}$. The
bounds on the singular values of $R^{\rho}_{ik}$, given up to a bounded factor
by the following lemma, give the bound
\[
\mu_\ell(\mathcal{L}_{ij}(s)\le Ce^{C|s|-\ell^{1/n}/C}\,,
\]
for some $ C$, which completes the proof of \eqref{eq:cest}.
\end{proof}

\begin{lem}
\label{lem:sing}
  Let $\rho\in(0,1)$ and $R_\rho:H^2(B_{\CC^n} ( 0, 1 ))\to H^2(B_{\CC^n} (
  0, \rho))$ induced by  the restriction map of $B_{\CC^n} ( 0, 1 )$ to
$B_{\CC^n} ( 0, \rho )$. Then, for any $\widetilde\rho\in(\rho,1)$ there exits a
constant
C such that
\[
\mu_\ell(R^\rho)\le C \widetilde\rho^{\ell^{1/n}}\,.
\]
\end{lem}
\begin{proof}
We use \eqref{eq:cec}
with the standard basis $(x_\alpha)_{\alpha\in\N^n}$ of $ H^2 ( B_{\CC^n} ( 0, 1 ) ) $: 
\begin{equation}
\label{eq:basis}
 x_\alpha (z) = c_\alpha z_1^{\alpha_1} \cdots z_n^{\alpha_n} \,, \ \ \
\int_{ B_{\CC^n} ( 0, 1 ) } |x_\alpha ( z ) |^2 dm ( z ) = 1 \,, \ \ 
\alpha \in {\N_0^n } \,, \end{equation}
for which we have
\[ \|R_\rho(x_\alpha)\|^2=\int_{ B_{\CC^n} ( 0, \rho ) } | x_\alpha ( w )|^2 dm ( w ) 
= \rho^{ 2(|\alpha| + n)} \,. \]
The number of $ \alpha$'s with $ |\alpha | \leq m $ is approximately
$ m^{n} $ and hence by \eqref{eq:cec} we have 
\[ \mu_\ell (  R_\rho) \leq C \sum_{C  k \geq \ell^{1/n}  }
k^{n-1} \rho^k  
\leq C \widetilde\rho^{\ell^{1/n}} \,.\]
\end{proof}
\medskip

The next proposition is a modification of standard 
zeta function arguments --  see \cite{Pol1} and \cite{Pol2} for the
discussion of the hyperbolic case.
\begin{prop}
\label{p:2}
Let $ {\mathcal L } ( s) $ be defined by \eqref{eq:trans}.
Then, if $ Z_\Gamma $ is the zeta function \eqref{eq:zeta} corresponding
to the group $ \Gamma $ then
\[ Z_\Gamma ( s) = \det ( I - {\mathcal L}( s ) ) \,. \]
\end{prop}
\begin{proof}

For $ s $ fixed and $ z \in \CC $
$$ h ( z ) \stackrel{\mathrm{def}}{=} \det ( I - z {\mathcal L} ( s ) ) $$
is, in view of \eqref{eq:cest} and \eqref{eq:weyl}, an entire function of
order $ 0$. For $ |z|$ sufficiently
small $ \log ( I - z {\mathcal L} ( s ) ) $ is well defined and we have
\begin{equation}
\label{eq:log} 
\begin{split}
\det ( I - z {\mathcal L} ( s ) ) & = \exp \left( - \sum_{ n=1}^\infty 
\frac{ z^n}{n} \tr ( {\mathcal L}(s))^n \right) \,.
\end{split}
\end{equation}

The correspondence between the closed geodesic (or, equivalently, 
conjugacy classes of hyperbolic elements) and the periodic orbits of 
$ T $  is particularly simple for
Schottky groups and we recall it in the form given in 
\cite{Pol2}
(where it is given in a more complicated setting of co-compact groups):
\begin{quote}
Closed geodesics on $ \Gamma \backslash \Hh $, $ \gamma $ of length
$ l( \gamma ) $, and word length $ |\gamma | $ are in one to one correspondence
with periodic orbits $ \{ x, Tx , \cdots , T^{n-1} x \} $ such that
 $  \| T^n ( x )\| = \exp \ell ( \gamma ) $, 
and $ n = |\gamma |$. 
For prime closed geodesics we have the same correspondence with primitive
periodic orbits of $ T $.
\end{quote}
It is not needed for us to recall the precise definition of the 
{\em word length}.
Roughly speaking it is the number of generators of $ \Gamma $ needed to 
write down $ \gamma $. 

To evaluate 
$ \tr {\mathcal L} ( s)^m $ we write 
\[  \tr  {\mathcal L} ( s)^m = 
\sum_{ (i_1 , \cdots , i_m ) , i_1 = i_m}
\tr \left( {\mathcal L}_{i_1 i_2 } \circ \cdots \circ {\mathcal L}
_{i_{m-1} i_m} \right) \,, \]
where in the notation of \eqref{eq:fij} we have 
\begin{gather*}
  {\mathcal L}_{i_1 i_2 } \circ \cdots \circ {\mathcal L}
_{i_{m-1} i_m}  u ( z ) = \| D ( f 
_{i_1 i_2 } \circ \cdots \circ {f}_{i_{m-1} i_m}  )(z) \|^s 
u ( f_{i_1 i_2 } \circ \cdots \circ {f}_{i_{m-1} i_m} (z) ) \,, \\ 
 f _{i_1 i_2 } \circ \cdots \circ {f}_{i_{m-1} i_m} \;: \; 
D_{i_1 } \ \longrightarrow \ D_{ i_1} \,. \end{gather*}
The trace of this operator is non-zero only if $ 
 f_{i_1 i_2 } \circ \cdots \circ {f}_{i_{m-1} i_m} $ has a fixed point
in $ D_{i_1} $. Since this transformation corresponds to an
element of $ \Gamma $ that fixed point is unique. Let us call this
element $ \gamma^{-1} $. Since it corresponds to 
a given periodic point, $ x $, of $ T^n $,
(corresponding to a fixed point of
$  f _{i_1 i_2 } \circ \cdots \circ {f}_{i_{m-1} i_m} $), $ \gamma $ 
is determined uniquely by $ x $ and $ n $:
\[ \gamma = \gamma ( x, n ) \,,  \ \ T^n x = x \,.\]
By conjugation and a choice of coordinates $ z $
it can be put into the form 
\[ \gamma ( z ) = e^{ \ell ( \gamma )} ( e^{i \theta_1 ( \gamma ) 
} z_1 , \cdots , e^{ i \theta_n ( \gamma )} z_n ) \,, \]
and the trace can be evaluated on the Hilbert space $ H^2 ( B_{ \CC^n}
( 0, 1 )) $. Using the basis \eqref{eq:basis} we can write the
kernel of $   {\mathcal L}_{i_1 i_2 } \circ \cdots \circ {\mathcal L}
_{i_{m-1} i_m} $ as 
\[  \begin{split}
 {\mathcal L}_{i_1 i_2 } \circ \cdots \circ {\mathcal L}
_{i_{m-1} i_m} ( z , w ) & = |\gamma' ( 0)|^{-s} 
\sum_{ \alpha \in \N_0^n } c_\alpha (\gamma^{-1}(z))^\alpha \bar w^\alpha 
\\ & = \sum_{ \alpha \in \N_0^n } c_\alpha 
e^{ - (s + |\alpha| )  \ell ( \gamma )  - i \langle \theta ( \gamma ) , 
\alpha \rangle } z^\alpha \bar w^\alpha \,. \end{split} \]
The evaluation of the trace is now clear.
 
Returning to \eqref{eq:log}, we obtain for $ \Re s $ sufficiently large
(using $ \{ \gamma \} $'s to denote the conjugacy classes of 
primitive of elements of $ \Gamma $), 
\[ 
\begin{split}
\det ( I - z {\mathcal L} ( s ) ) & = \exp \left( - \sum_{ n=1}^\infty 
\frac{ z^n}{n}  \sum_{ T^n x = x  }   
 \sum_{ \alpha \in \N_0^n } 
e^{ - (s + |\alpha| )  \ell ( \gamma (x,n) )  - i \langle \theta ( 
\gamma (x,n)) , 
\alpha \rangle }  
\right)  \\
& = 
 \exp \left( - \sum_{ n=1}^\infty 
  \sum_{ { \{\gamma\} }\atop{|\gamma| = n}}  \sum_{k=1}^\infty 
\frac{z^{nk}} {k}  
 \sum_{ \alpha \in \N_0^n } e^{ - k( 
(s + |\alpha| )  \ell ( \gamma )  - i \langle \theta ( \gamma ) , 
\alpha \rangle ) }  
 \right) \\
& =  \prod_{ {\{ \gamma \}}  } \prod_{\alpha \in \N_0^n} 
\left( 1 -  z^{|\gamma|} e^{ - i \langle \theta( \gamma ) , \alpha \rangle }
e^{ - ( s + |\alpha| )  \ell( \gamma ) } \right)
\end{split}
\]
which in view of \eqref{eq:zeta} proves the proposition once
we put $ z = 1 $. 
\end{proof}

\medskip

\noindent
{\bf Remark.} The proof above is inspired by the work on
the distribution of resonances in Euclidean scattering - see
\cite[Proposition 2]{Zw}. The
Fredholm determinant method and the use of Weyl inequalities in the
study of resonances
were introduced by Melrose \cite{Mel} and developed further by many authors
-- see \cite{Sj},\cite{MZ1}, and references given there. That was
done at about the same time as David Fried (across the Charles
River from Melrose) was applying the Grothendieck-Fredholm theory
to multidimensional zeta-functions \cite{fr}.
In both situation the enemy is the exponential
growth for complex energies $ s $, which is eliminated thanks to
analyticity properties of the kernel of the operator.

Finally, we remark that 
in view of the lower bounds on the number of zeros of $ Z _\gamma $
obtained in \cite{gz3} in dimension two, and in 
\cite{Pe} in general, we see from Proposition \ref{p:2} that
the upper bound \eqref{eq:det} is optimal for any $ \gamma $.

\section{Quasi-self-similarity of $ \Lambda ( \Gamma ) $ 
and Markov partitions for $ \Gamma \backslash \Hh^{n+1} $} 
\label{mp}

\begin{figure}
\begin{center}
\includegraphics[width=3in]{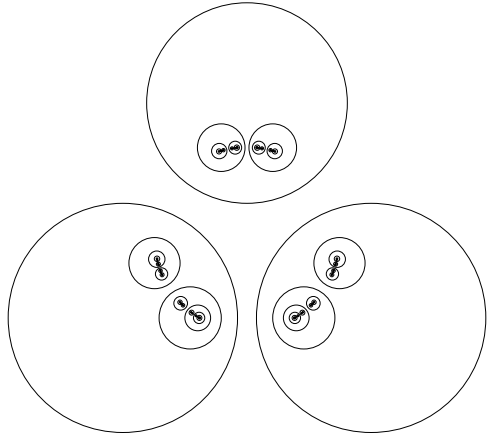}
\end{center}
\caption{A Markov partition and its refinements for  
the reflection group $\Gamma_\theta$ of Fig. 1 (considered here as a Kleinian
conformal group  on $\partial\Hh^3\simeq\RR^2\cup\infty$).}
\label{fi:1}
\end{figure}

In this section we will review the results on convex co-compact Schottky
groups coming essentially from \cite{Sull} and \cite{McM}.  The geometric
point of view presented here was explained to the authors by Curt
McMullen\footnote{who also generously provided us with Figures
  \ref{fi:2},\ref{fi:1}, and previously with the essential part of Figure
  \ref{fi:0}.}.

We start with a more general definition of convex co-compact subgroups
of $ {\mathrm{Isom}} ( \Hh^{n+1} ) $. A discrete subgroup is called
{\em convex co-compact} if 
\begin{equation}
\label{eq:coco}   \Gamma \backslash C( \Lambda ( \Gamma )) 
\ \ \text{is compact} \,, \ \  
C ( \Lambda ( \Gamma ) ) \stackrel{\mathrm{def}}{=}  
 \text{convex hull}  ( \Lambda ( \Gamma ) )\end{equation}
Here, the convex hull is meant in the sense of the hyperbolic 
metric on $ \Hh^{n+1} $: $ \Lambda ( \Gamma ) \subset \partial \Hh^{n+1} $,
and $ \gamma $ acts on it in the usual way. In particular this implies
that $ \Gamma \backslash
 C ( \Lambda( \Gamma ) ) $ has a compact fundamental domain in $ \Hh^{n+1} $. 

The first result gives a quasi-self-similarity for arbitrary 
convex co-compact groups:

\begin{prop}
\label{p:3}
Suppose that $ \Gamma \subset {\mathrm{Isom}} ( \Hh^{n+1} ) $ is 
convex co-compact in the sense of \eqref{eq:coco}.
Then there exist $ c > 0 $ and $ r_0 > 0 $ 
such that for any $ x_0 \in \Lambda ( \Gamma ) $ 
and $ r < r_0 $ there exists a map $ g: B_{\SP^n }(x_0, r)
 \rightarrow \SP^n $ with the properties
\begin{gather}
\label{eq:ss}
\begin{gathered}  g ( \Lambda ( \Gamma ) \cap  B_{\SP^n}( x_0, r)  ) 
\subset \Lambda( \Gamma ) \\
 c r^{-1} d_{\SP^n} (x,y)  \leq d_{\SP^n}( g ( x ) ,g ( y ) )
 \leq c^{-1} r^{-1}  d_{\SP^n} (x,y)
\,, \ \ x, y \in B_{\SP^n} ( x_0, r) \,.
\end{gathered}
\end{gather}
\end{prop}
\begin{proof} 
We proceed following the argument  
in \cite[Sect.3]{Sull}. 
Let us fix $ z_0 \in  C ( \Lambda ( \Gamma )) $.  If $ L $ is the geodesic ray through $ z_0 $
and $ x_0 $, then 
\[ \exists \; C> 0 \; \forall z \in L \; \exists \;  \gamma \in \Gamma 
\ \ d ( \gamma^{-1}  z_0 , z ) < C \,.\]
This follows from the compactness of $  \Gamma  \backslash
C ( \Lambda ( \Gamma ) )$:
for any point on the ray, $ z $, 
there exists an element of the orbit of $ z_0 $
within a finite distance from $ z $. We can now choose 
$ z = z ( r ) 
 $ on the ray $L$  so that $ d ( z, z_0 ) = \log ( 1 / r )$, 
and then $ \gamma $ such that $ d ( \gamma^{-1} z_0 , z_0 ) = 
\log ( 1 / r ) + {\mathcal O } ( 1 ) $. 

If $ x_\gamma $ is the end point of the geodesic ray through $ z_0 $ and 
$ \gamma^{-1} z_0 $, then for a fixed $ C_1 $, the ball 
$ B_{\SP^n } ( x_\gamma, C_1 r ) $ covers $ B_{ \SP^n } ( x_0 ,r ) $.
The action of $ \gamma $ 
on $ B_{\SP^n } ( x_\gamma, C_1 r ) $ satisfies \eqref{eq:ss}: 
$ \Lambda ( \Gamma ) $ is $ \Gamma$-invariant, and the other property 
follows by putting $ \gamma $ into the normal form \eqref{eq:hyp}.
Since $ z_0 $ was fixed and we have no dependence on $ \gamma $,
the proof is completed.
\end{proof}

Strictly speaking we will not use the quasi-self-similarity
explicitely. It is present implicitely in the proof of Proposition 
\ref{p:4} below where again the compactness of $ \Gamma \backslash
C \Lambda ( \Gamma ) $ is crucial. We remark that for $ n = 1 $,
that is for surfaces, Proposition \ref{p:3} shows that, in the
upper half plane model chosen so that $ \Lambda \Subset \RR $,
\[ \Lambda ( \Gamma ) + [ -h, h] = \bigcup_{ j=1}^{N ( h ) } 
[ a_j (h ) , b_j ( h ) ] \,, \ b_j ( h ) < a_{ j+1} ( h ) \,, \ \ 
|b_j ( h ) - a_j ( h )| \leq K h \,,\]
with $ K $ independent of $h $. That gives a simple proof 
of our main result -- see Section \ref{etd}.

For the general 
we recall the definition of a Markov partition in a
simple form \cite{McM} applicable here. Thus let $ \mu_\Lambda $ 
be the $ \delta$-Hausdorff measure restricted to $ \Lambda ( \Gamma )$. 
Let $ T $ and $ \mathcal{D} $ be defined as in \eqref{eq:T}.

A {\em Markov partition} for the  Schottky
group $\Gamma$  is given by 
a collection, $ {\mathcal P} $, of subsets of $ \mathcal{D} \subset \SP^n $,
satisfying
\begin{gather}
\label{eq:markov}
\begin{gathered}
\forall \; P,Q  \in {\mathcal P} \ \ \ 
 \ { \mu_\Lambda ( T( P ) \cap Q ) \neq 0 }  
\ \Longrightarrow \ T( P ) \supset   Q 
\\
\text{$T$ is a homeomorphism on a neighbourhood of $ P \cap 
(T\rest_P ) ^{-1} ( Q ) $ if $ \mu_\Lambda ( T ( P ) \cap Q ) \neq 0 $} \\
\forall \; P \in {\mathcal P} \ \ \ \mu_\Lambda ( P ) > 0 \\
\forall \; P,Q \in {\mathcal P}\,, \ P \neq Q  \ \ \ \ \ 
\mu_\Lambda ( P \cap Q ) = 0   \\
\mu_\Lambda ( T ( P ) ) =  \sum_
{{Q \in {\mathcal P} }\atop{ \mu_\Lambda ( T( P ) \cap Q ) \neq 0}  }
\mu_\Lambda ( Q ) \,.
\end{gathered}
\end{gather}

A refinement of $ {\mathcal P} $ is a new Markov partition given by 
\[ {\mathcal R} ( {\mathcal P} ) \stackrel{\mathrm{def}}{=} 
\{  (T\rest_P ) ^{-1} ( Q ) \;: \; P, Q \in {\mathcal P} \,, \ \
\mu_\Lambda ( T(P) \cap  Q ) \neq 0 \} \,, \]
see Fig.\ref{fi:1}. The refinements can of course be iterated giving
$ {\mathcal R}^m ( {\mathcal P} ) $. 

The Markov partition 
is said to be expanding if there exists a metric $\|\hphantom{v}\|$ 
on $\mathcal{D}$ and a real $\theta>1$
such that $\|DT(v)\xi\|\ge\theta\|\xi\|,v\in V,\xi\in T_vV$.
An example of a Markov partition  for the Schottky group map is given by 
\begin{equation}
\label{eq:mare}
 {\mathcal P_0} 
 = \{ \mathcal{D}_l \;: \; 1 \leq i \leq 2\ell \} \,. 
\end{equation}
If $N$ is big enough its refinement, $\mathcal{D}^N = 
{\mathcal R}^N ( {\mathcal P} ) $, 
 considered in Section 4, is expanding.

The result we need to estimate the zeta function is given in the following
\begin{prop}
\label{p:4}
Let $ \Gamma $ be a convex co-compact Schottky group as defined in 
Section 2. Then there exist positive constants $ h_0 $, $ C_0 $, and 
$ C_1 $, 
such that for $ 0 < h < h_0 $ 
there exists a Markov partition, $ {\mathcal P} ( h ) $,
for $ \Gamma $ with the following
properties:
\begin{gather}
\label{eq:mar2}
\begin{gathered}
\forall \; P \in {\mathcal P} (h) \; \exists \; 
x_P \in \SP^n \,, \    \ \ \ \
P \subset B_{\SP^n }  ( x_P , C_0 h  ) \\
\forall \; P,Q\in {\mathcal P} (h) \ \ \ \ \Lambda ( \Gamma ) 
\cap  T ( P ) 
\cap Q  \neq \emptyset  \ \Longrightarrow \ 
d_{ \SP^n } (( T \rest_P)^{-1} ( Q) , \partial P ) > h / C_1 \,.
\end{gathered}
\end{gather}
In fact, we can make each element of $ \mathcal{P}(h) $ a union of a fixed number
of balls of radii comparable to $ h$.
\end{prop}
\begin{proof}
Rather than use the successive refinements of the Markov partition 
$\mathcal{D}^N$ (which would not work) 
we will apply the following geometric observation:
the projection of any geodesic in $ C( \Lambda ( \Gamma ) ) $ to $ \Gamma 
\backslash \Hh^{n+1} $ intersects the projections of $ C ( \partial \mathcal{D}_l ) $ 
at uniformly bounded intervals, say by $ d_0 > 0$. This follows from the
bound on the diameter of the compact fundamental domain of
$ \Gamma \backslash C ( \Lambda ( \Gamma ) ) $.

For $ e^{-m} < h < e^{-m+1} $ we choose ${\mathcal P  } (h) $ in three 
steps. Let $ z_0 \in \Hh^{n+1} $ be a fixed point, and $ {\mathcal P}_0 $ the
partition $\mathcal{D}^N$. We then put 
\begin{gather*}
 {\mathcal P}_k (h) =
 \{ P  \in {\mathcal R}^k
( {\mathcal P}_0 ) \;: \;  
  m < d_{\Hh^{n+1}} (C( P) , z_0 ) \leq m + d_0  \} 
 \\
\widetilde {\mathcal P} (h) \stackrel{\mathrm{def}}{=}  
\bigcup_{ k\geq 1 } 
\{  P \in  {\mathcal P}_k (h) \;: \; P \not \subset Q \in 
 {\mathcal P}_{k-1} (h) \} \,.
\end{gather*}
The choice of $ d_0 $ shows that the elements of $ \widetilde {\mathcal P}
(h ) $ cover $ \Lambda ( \Gamma ) $. From  $ d ( C ( P ) , z_0 ) ) = 
m + {\mathcal O} ( 1 ) $, we conclude that 
$ P \subset B_{\SP^n} ( x_P , C_0 h ) $, $ C_0 $ uniformly bounded, 
$ h \sim e^{-m } $.

To obtain the strong form of the Markov partition
property given in the second part of \eqref{eq:mar2} we modify 
$ \widetilde {\mathcal P} (h) $ as follows:
\begin{equation}
\label{mard}
\widehat{{\mathcal P}} ( h ) = 
\big\{ \{  x \in \SP^n \;: \; d( x , \Lambda ( \Gamma ) \cap P ) < h/C_1 \} \;: \;
P \in \widetilde{\mathcal P} ( h ) \big\}\,.
\end{equation}
We still have the property that
\[ \forall \; P, Q \in \widehat{{\mathcal P}} ( h )  \ \ 
P \ne Q \ \Longrightarrow \ 
P \cap Q \cap \Lambda ( \Gamma ) = \emptyset \,.\]
If $ \widehat{{\mathcal P}} (h) \ni Q  \subset D_j $, 
where $ D_j $'s are as \eqref{eq:fij}
then $ f_{ij} $ is a contraction of $ Q $ into a 
$ \widehat{{\mathcal P}} (h) \ni Q'  \subset D_i $.
In fact, $ f_{ij} $ is a contraction
of $ D_j $ into $ D_i $ and it preserves $ \Lambda ( \Gamma ) $. 
That shows the second part of \eqref{eq:mar2}.

Depending on the contraction constants of $ f_{ij}$'s (a finite set) 
we can replace elements of $ \widehat {\mathcal P}  (h ) $ with unions of
a fixed number of balls (the number depending only on the constants $ C_0$
and $ C_1 $ and the dimension) so that all the properties still hold.
\end{proof}

\section{Estimates in terms of the dimension of $ \Lambda ( \Gamma ) $.}
\label{etd}

In Proposition \ref{p:2} we used the Markov partition 
\eqref{eq:mare}. It is clear from the proof that a Markov partition 
for which \eqref{eq:mar2} holds (for some fixed $ h$) can be used
instead as long as the contraction property is valid in some complex 
neighbourhood of the components of the partition. 
For the proof of Theorem stated in Section 1 we will modify $ D_j $'s in
the definition of $ {\mathcal L } ( s)  $ in a way dependent on the
size of $ s $: we will use Proposition \eqref{p:4} with  $ h = 1/ |s| $. 
The self-similarity structure of 
$ \Lambda ( \Gamma ) $ will show that we can choose $ D_j = D_j ( h ) $
to be a union of $ {\mathcal O} ( h^{-\delta } ) $ disjoint balls of 
radii $ \sim h $.
A modification of 
the argument used in the proof of Proposition \ref{p:1} will then 
give \eqref{eq:zetad}. 

We start with the following lemma which is a more precise version of
the argument already used in the proof of Proposition \ref{p:1}:
\begin{lem}
\label{l:5}
Suppose that $ \Omega_j \subset \CC^n $,
$j=1,2$, are open sets, and
$ \Omega_1 = \bigcup_{k=1}^K B_{\CC^n } ( z_k, r_k) $. Let 
$ g $ be a holomorphic mapping defined on a neighbourhood, $ \widetilde 
\Omega_1 $ of
$ \Omega_1 $ 
with values in $\Omega_2$, satisfying 
\[  d_{ \CC^n } ( g ( \Omega_1 ) , \partial \Omega_2 ) > 
{1}/{C_0} 
> 0 \,.\]
If 
\[  A \;: \; H^2 ( \Omega_2 ) \ \longrightarrow H^2 ( \Omega_1 ) 
\,, \ \ A u ( z ) \stackrel{\mathrm{def}}{=} u ( g ( z ) ) \,, \]
then for some $ C_1 $ depending only on $ r_k$'s, $ K$, 
$ d_{\CC^n } ( \Omega_1 , \partial \Omega_2 )$, and $ \min_{\widetilde 
\Omega_1 }  \|Dg \|_{ \CC^n \rightarrow \CC^n} $, we have
\[ \mu_\ell (  A ) \leq  C_1 e^{ - \ell^{1/n} / C_1 } \,, \]
where $ \mu_{\ell} (A)$'s are the characteristic values of $ A $.
\end{lem}
\begin{proof}
We define a new Hilbert space 
\[ \HH \stackrel{\mathrm{def}}{=} \bigoplus_{k=1}^K H^2 ( B_k ) \,, \ \ 
B_k = B_{ \CC^n } ( z_k , r_k)  \,, \]
and a natural operator
\[ J \;: \; H^2 ( \Omega_1 ) \; \longrightarrow \; {\mathcal H} \,,
\ \ ( Ju)_k  = u\rest_{B_k} \,.\]
We easily check that $ J^* J: H^2 ( \Omega_1 ) \rightarrow H^2 ( \Omega_1 ) $
is invertible with constants depending only on  $K$. Hence 
\[ \mu_\ell ( A ) = \mu_\ell ( ( J^* J)^{-1} J^* J A ) \leq 
\|  ( J^* J)^{-1} \|\, \|J^* \| \mu_\ell ( J A ) \,. \]
We then notice that 
\[ \mu_{k\ell} ( J A ) \leq k \max_{ 1 \leq k \leq K } \mu_\ell ( A_k ) \,,\]
where 
\[ A_k \;: \; H^2 ( \Omega_2 ) \; \longrightarrow \; H^2 ( B_k ) \,, 
\ \  A_k u ( z) = u ( g_k ( z) ) \,, \ \ g_k = g\rest_{ B_k } \,. \]
To estimate the characteristic values of $ A_k $ we observe that
we can extend $ g_k $ to a larger ball, $ \widetilde B_k $
(contained in $ \widetilde \Omega_1 $)
and such that the image of its closure still lies in $ \Omega_2 $ (since
we know that $ \min_{ \widetilde \Omega_1 }
 \|Dg \|_{ \CC^n \rightarrow \CC^n}  $ is strictly less than
$1 $). That gives us the operators $ R_k: H^2 ( \widetilde B_k ) \rightarrow
H^2 ( B_k ) $, $ R_k u  = u\rest_{B_k } $, and $ \widetilde A_k$
defined as $ A_k $ but with $ B_k $ replaced by $ \widetilde B_k $. 
We now have $ A_k = R_k \widetilde A_k $ and consequently,
\[ \mu_\ell ( A_k ) \leq \| \widetilde A_k \| \mu_\ell ( R_k ) \,.\]
Lemma \ref{lem:sing} gives  $ \mu_\ell ( R_k ) \leq C_2 \exp ( -\ell^{1/n} /C_2 ) $
completing the proof.
\end{proof}

\medskip
\noindent
{\em Proof of Theorem.}
As outlined in the beginning of the section we put $ h = 1/|s| $, where
$ |s| $ is large but $ | \Re s | $ is uniformly bounded. In Proposition
\ref{p:4} each element of the Markov partition is given as union of
(a fixed number of) balls. We can complexify this set by taking a 
corresponding union of balls in $ \CC^n $ and all the properties hold
for the analytic continuations of $ f_{ij}$'s defined in \eqref{eq:fij}. 

The now classical results of Patterson and Sullivan \cite{Sull}
on the dimension of
the limit set show that the total number of the balls 
is $  {\mathcal O} ( h^{-\delta })$: what we
are using here is the fact that the Hausdorff measure of $ \Lambda ( \Gamma ) $
is finite. 

We can now apply the same procedure as in the proof of Proposition \ref{p:1}
using Lemma \ref{l:5}.
What we have gained is a bound on the weight: since $ |\Re s | \leq C$ and
$ f_{ij}' $ is real 
on the real $\SP^n$
\[ | [f_{ij }' (z) ]^s | \leq C \exp (  |s| |\arg  f_{ij }' (z)| ) \leq
C \exp ( C_1 |s| |\Im z | )  \leq C_2 \,, \ \ z \in D_j ( h ) \,. \] 
We write $ {\mathcal L} ( s) $ as a sum of fixed number of 
operators $ {\mathcal 
L}_{ij} ( s ) $ each of which is a direct sum of 
$ {\mathcal O} ( h^{-\delta} )  $ operators. 
The 
balls and contractions are uniform after rescaling by $ h$ and hence
the characteristic values of each of these operators satisfy the bound
$ \mu_\ell \leq C \gamma^l $, $ 0 < \gamma < 1 $. 
Using \eqref{eq:weyl} and \eqref{eq:obv} we obtain the bound
\[ \log | \det ( I - {\mathcal L} ( s) ) |  \leq C P(h)  = {\mathcal O} ( 
h^{-\delta}) \,, \]
and this is \eqref{eq:zetad}.
\stopthm
  
\medskip
\noindent
{\em Proof of Corollary 1.}
The definition of $ Z_\Gamma ( s) $ \eqref{eq:zeta} shows that 
for $ \Re s > C_1 $ we have $ |Z_\Gamma ( s ) | > 1/2 $. The Jensen formula
then shows that the left hand side of \eqref{eq:ubd} is bounded by
\[  \sum \{ 
m_\Gamma (s) \;: \;  | s - ir -  C_1 | \leq C_2 \} \leq 
2 \max_{ {|s| \leq r + C_3 }\atop{|\Re s |\leq C_0 }} \log |Z_\Gamma ( s ) |
+ C_4 \,,\]
and \eqref{eq:ubd} follows from \eqref{eq:zetad}.
\stopthm

\section{Schottky manifolds and resonances}

We recall that s complete Riemannian manifold of
constant curvature $-1$ is said to be Schottky if its
fundamental group is Schottky. In low dimensions Schottky manifolds can be
described geometrically. 

\begin{prop}
\label{prop:ScFu}
 Any convex co-compact hyperbolic surface is Schottky.
\end{prop}
This result is proved by Button \cite{Bu} and for the 
reader's convenience we
sketch the proof.
\begin{proof}
  Any convex co-compact surface is topologically described by two integers
  $(g,f)$~: its numbers $g$ of holes and $f$ of funnels, with the
  conditions $g\ge0$, $f\ge1$ and $f\ge3$ if $g=0$. For any such pair
  $(g,f)$, there does exist a Schottky surface of this type
 and we choose for each type $(g,f)$ such a
  surface $M_{g,f}$. The projection onto $M_{g,f}$ of the boundary of the
  Schottky domain 
is a collection
  $\mathcal{L}_1,\ldots,\mathcal{L}_\ell$, of mutually disjoint geodesic
  lines.
  
Let $M$ be any hyperbolic convex co-compact surface . The surface $M$ is
homeomorphic to some $M_{g,f}$. Pushing back on $M$ the geodesic lines
$\mathcal{L}_i,i=1,\ldots,\ell$ of $M_{g,f}$ and cutting $M$ along these
curves, we obtain in the hyperbolic plane a domain whose boundary is the
union of paired mutually disjoint curves
$\mathcal{C}_i,\mathcal{C}_{\ell+i},i=1,\ldots,\ell$, each one with a pair
of points at infinity. These points pair determine intervals, which are
mutually disjoint (the curves $\mathcal{C}_j,j=1,\ldots,2\ell$ don't intersect). The
intervals are paired with an hyperbolic transformation, so give a Schottky
group, which coincide with the fundamental group of the surface $M$.
\end{proof}

\medskip
\noindent
{\em Proof of Corollary 2.}
For a Schottky manifold, the fundamental groups is Schottky, and hence,
$M=\Gamma\backslash\Hh^{n+1},\Gamma\subset\mathrm{Isom}^+(\Hh^{n+1})$. 
We then 
introduce its zeta function $Z_M$ as the zeta function $Z_\Gamma$ of the
group $\Gamma$. Following Patterson and Perry \cite{PP} we introduce the
spectral sets $\mathcal{P}_M$ and $\mathcal{S}_M$ defined by the
Laplace-Beltrami operator $\Delta_M$ on $M$:
\begin{eqnarray*}
  \label{eq:Zsp}
  \mathcal{P}_M&=&\{s \; : \; 
\Re s>n/2, s(n-s)\hbox{ is a $L^2$ eigenvalue of $\Delta_M$}\}\\
  \mathcal{S}_M&=&\{s \; : \; \Re s<n/2, s\hbox{ is a singularity of the scattering matrix ${S}_M$ }\}
\end{eqnarray*}
Moreover, each complex $s$ in $\mathcal{P}_M$  has a multiplicity denoted by
$m_M(s)$, each $s$ in $\mathcal{S}_M$ a pole multiplicity denoted by
$m^-_M(s)$. 
In the case of surfaces,
the 
divisor of the Selberg zeta function $Z_M$ is given by the following formula:
\[
-\chi_M \sum_{k=0}^\infty(2k+1)[-k] +
m_M\left(\frac{n}2\right)\left[\frac{n}2\right]
+\sum_{s\in \mathcal{P}_M}m_M(s)[s]+\sum_{s\in  \mathcal{S}_M}m^-_M(s)[s]\,,
\]
where $\chi_M$ is the Euler-characterictic of $M$, see 
\cite[Theorem 1.2]{PP}  The zeta function, $Z_M$,
is entire and in any half-plane $\{\Re s>-C_0\}$, and 
the formula above shows that the bounds on the number of its zeros
provide bounds on the number of resonances. The dimension of 
the limit set, $ \delta $ depends only on $ \Gamma $ and, as shown in 
\cite{Sull}, it gives the Hausdoroff dimension of the 
recurrent set for the geodesic flow on $T^*M$ by the formula $ 2 ( 1 + 
\delta ) $. 
\stopthm

For a convex co-compact hyperbolic manifold $M$, Patterson and Perry give a
formula for the divisor of the zeta function $Z_M$ in any (even) 
dimension, but it
does not imply (in the non-Schottky case) 
that the zeta function is entire. In the case of Schottky
groups, the zeta function $Z_M$ is entire, as it was shown in the Proposition
\ref{p:2}. Hence 
we concluded that Corollary 2 holds also for Schottky manifold.

\medskip

We conclude with some remarks about Kleinian 
groups in dimension $n+1=3$. 
Schottky 3-manifold are geometrically described by Maskit
\cite{Mas1}:
\begin{prop}
  A hyperbolic convex co-compact, non compact 3-manifold is Schottky if
  and only if its fundamental group is a free group of finite type.
\end{prop}

While non compact surfaces of finite geometric type have always a free
fundamental group, that is not true for 3-manifold. For instance,
if $\Gamma$ is a
co-compact surface group, the 3-manifold $\Hh^3/\Gamma$ is convex co-compact
with a non-free fundamental group. Quasi-fuchsian groups (that is,
deformation of
such a $\Gamma$ in $\mathrm{Isom}(\Hh^3)$) give similar examples.

Finally, we note  that the bound on the number of zeros of $Z_\Gamma$
 established here 
for Schottky groups is valid for any group $\Gamma$, for which a expanding
Markov partition can be built.  Anderson and Rocha \cite{AR} construct such
Markov partition for any function group.  This class of groups does not
exhaust all convex co-compact groups (the complement in the 3-sphere of a
regular neighbourhood of a graph is not in this class) and it is not known
if all convex co-compact Kleinian groups admit an expanding Markov partition.

\section{Numerical Results.}
\label{num}

\subsection{Discussion.}

In this section, we present numerical results on the distribution of
zeros of a closely related dynamical zeta function $Z(s)$ for a
conformal dynamical system. We consider the simple case of groups $
\Gamma_\theta $ generated by reflections in three symmetrically placed
circles perpendicular to the unit circle (considered as the boundary of
the Poincar\'e disc), and intersecting it at angles $ \theta $ -- see
Fig.\ref{fi:0} where $ \theta = 110^\circ$.

Numerical computations of the zeta function in that case have been
already performed by McMullen \cite{McM} and Jenkinson-Pollicott
\cite{Jen}. Their goal was to find an efficient way of computing the
dimension of limit sets.  Table {\ref{tab:dims}} gives the (approximate)
dimensions of the limit sets for the relevant angles, calculated as the
largest real zero of $Z(s)$ {\cite{Jen}} using Newton's method.

\begin{figure}
\begin{center}
\PSbox{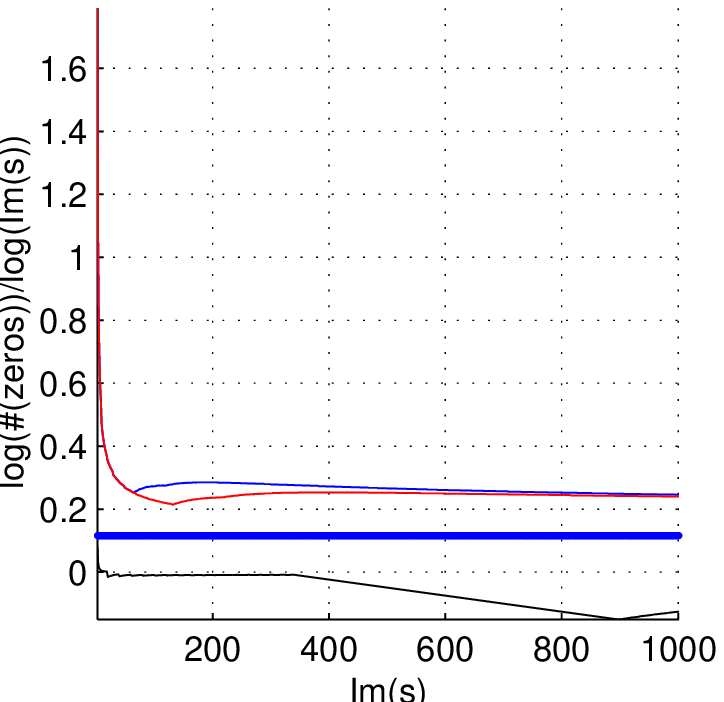}
\end{center}

\caption{This plot shows $\frac{\log(|\{s\in[x_0,x_1]\times[y_0,y] :
Z(s)=0\}|)} {\log(y)} - 1$ as a function of $y$, for different values of
$x_0$: The thin blue line is for $x_0=-0.2$, the red line for
$x_0=-0.1$, and the black line for $x_0=+0.1$.  Note that the value of
$x_1$ is not very important because $Z(s)$ decays very rapidly for large
$\Re(s)$.  Thus, we set $x_1=10$ throughout.  The value of $y_0$ is fixed
at $-0.1$, to avoid integrating over any zeros.  The thick horizontal
line indicates the dimension of the corresponding limit set.  In this
plot, $\theta=10^\circ$.}

\label{fig:run1-10}
\end{figure}

\begin{figure}
\begin{center}
\PSbox{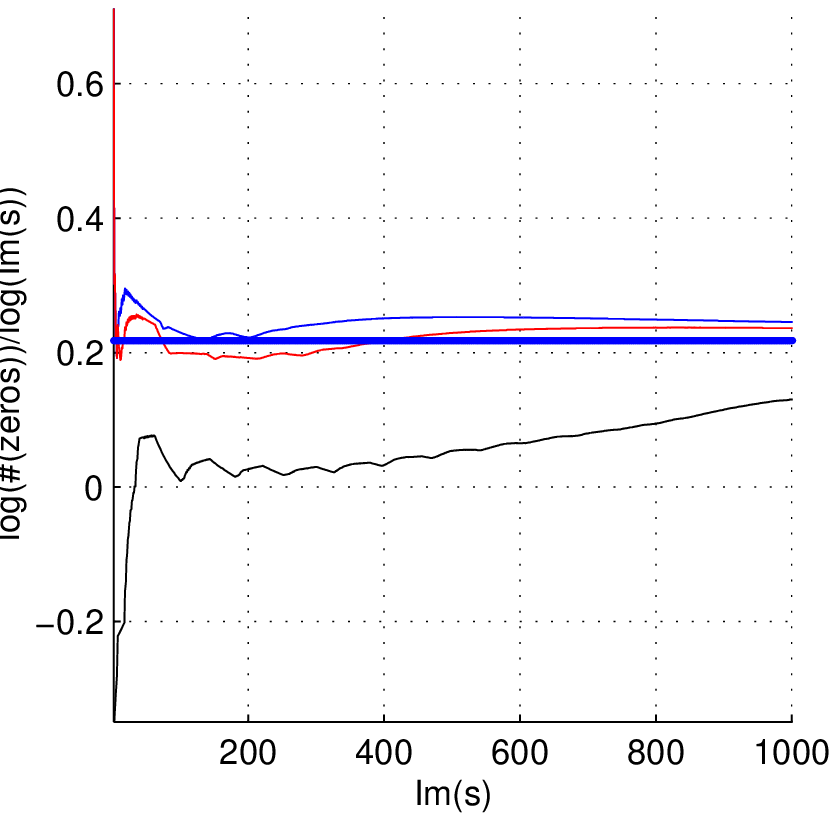}
\end{center}

\caption{This plot shows $\frac{\log(|\{s\in[x_0,x_1]\times[y_0,y] :
Z(s)=0\}|)} {\log(y)} - 1$ as a function of $y$, for different values of
$x_0$: The thin blue line is for $x_0=-0.2$, the red line for
$x_0=-0.1$, and the black line for $x_0=+0.1$.  Note that the value of
$x_1$ is not very important because $Z(s)$ decays very rapidly for large
$\Re(s)$.  Thus, we set $x_1=10$ throughout.  The value of $y_0$ is fixed
at $-0.1$, to avoid integrating over any zeros.  The thick horizontal
line indicates the dimension of the corresponding limit set.  In this
plot, $\theta=40^\circ$.}

\label{fig:run1-40}
\end{figure}

Figures {\ref{fig:run1-10}}-{\ref{fig:run1-40}} show
\begin{equation}
  \frac{\log(|\{s\in[x_0,x_1]\times[y_0,y] : Z(s)=0\}|)} {\log(y)} - 1
\end{equation}
as a function of $y$.  In each plot, the value of $x_0$ is varied to
test the dependence of the distribution on the region in which we count:
The blue line corresponds to $x_0=-0.2$, the red line $x_0=-0.1$, and
the black line $x_0=+0.1$.  The data show that most of the zeros lie in
the left half plane.  Based on the theorems proved in earlier sections,
we expect the curves to be bounded above by the dimension (the thick
blue line) asymptotically.  This is not the case, except for the black
line, which represents zeros with $\Re(s)>x_0=+0.1$.

\begin{figure}
\begin{center}
\PSbox{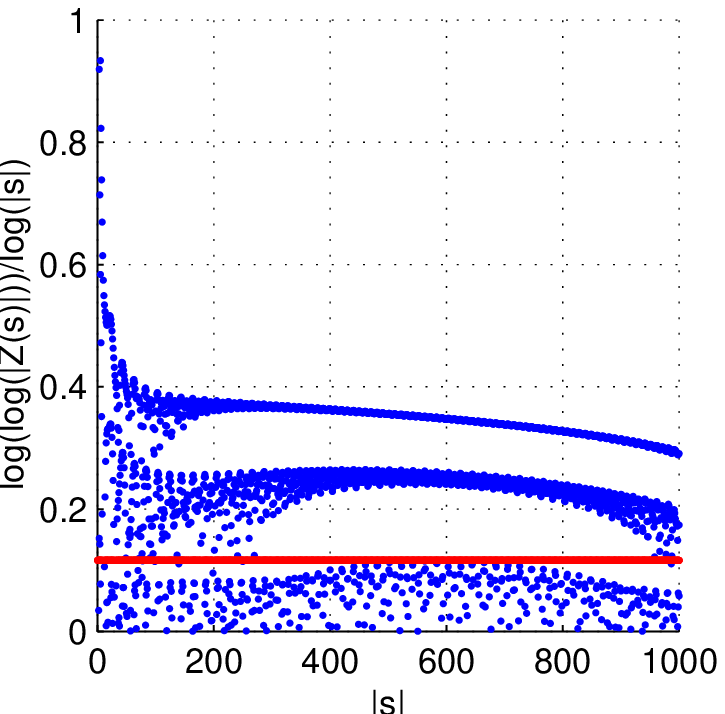}
\end{center}

\caption{This plot shows $\log(\log(|Z(s)|))/\log(|s|)$ for a large
number of points in the rectangle $[-0.2,1]\times[0,10^3]$.  The
horizontal line indicates dimension.  Here, $\theta=10^\circ$.}

\label{fig:zrad-10}
\end{figure}

\begin{figure}
\begin{center}
\PSbox{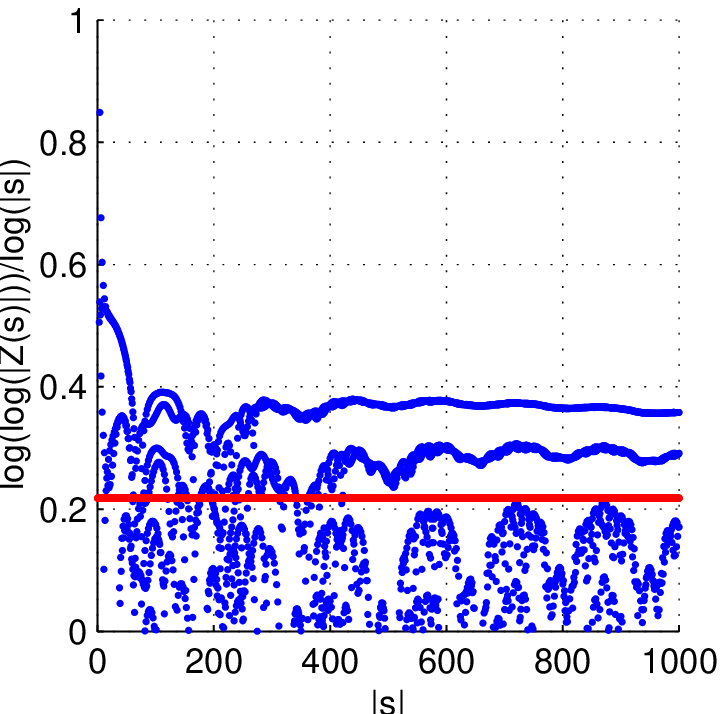}
\end{center}

\caption{This plot shows $\log(\log(|Z(s)|))/\log(|s|)$ for a large
number of points in the rectangle $[-0.2,1]\times[0,10^3]$.  The
horizontal line indicates dimension.  Here, $\theta=40^\circ$.}

\label{fig:zrad-40}
\end{figure}

Similarly, Figures {\ref{fig:zrad-10}}-{\ref{fig:zrad-40}} show
$\frac{\log(\log(|Z(s)|))}{\log(|s|)}$ as a function of $|s|$, for a
large number of points in the rectangle $[-0.2,1.0]\times[0,10^3]$.  In
this case, we also expect the curves to be asymptotically bounded by the
dimension.  This is also not the case.  The only reasonable explanation,
barring errors in the numerical calculations, is that the asymptotic
upper bound is accurate only for very large values of $\Im(s)$, and we
were not able to calculate $Z(s)$ reliably for such values.  These
results also show that $Z(s)$ has plenty of zeros in regions of
interest.

\begin{table}

\begin{center}
\begin{tabular}{c|c}
$\theta$ & dim\\\hline
$10^\circ$ & 0.11600945\\
$20^\circ$ & 0.15118368\\
$30^\circ$ & 0.18398306\\
$40^\circ$ & 0.21776581\\
\end{tabular}
\end{center}

\caption{Dimensions of the limit set for relevant values of $\theta$.}
\label{tab:dims}
\end{table}

\begin{figure}
\begin{center}
\PSbox{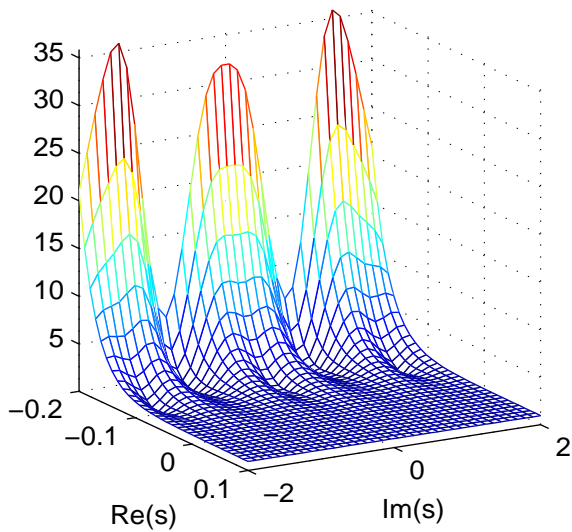}
\end{center}

\caption{This plot shows $|Z(s)|$ in the square
$[-0.2,0.3]\times[-2,2]$, for $\theta=30^\circ$.}

\end{figure}

\begin{figure}
\begin{center}
\PSbox{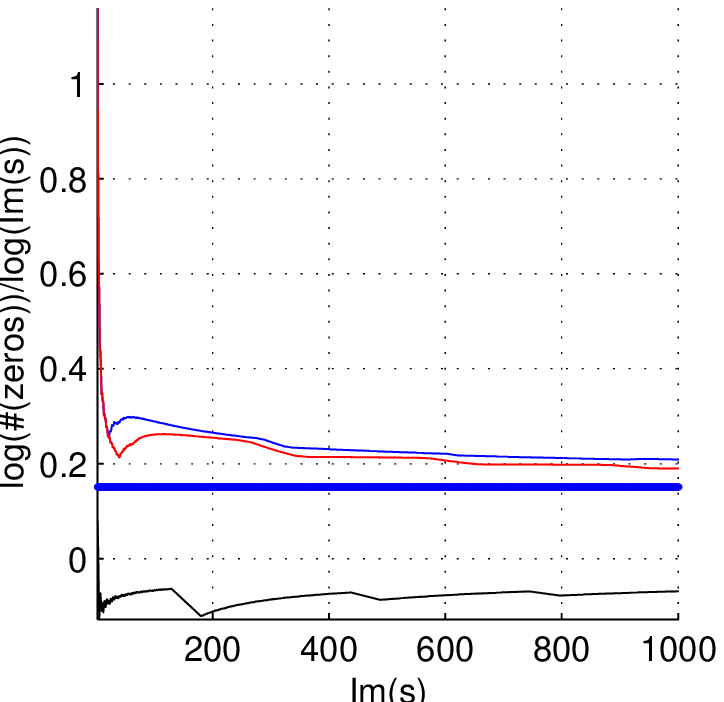}
\end{center}

\caption{This plot shows $\frac{\log(|\{s\in[x_0,x_1]\times[y_0,y] :
Z(s)=0\}|)} {\log(y)} - 1$ as a function of $y$, for different values of
$x_0$: The thin blue line is for $x_0=-0.2$, the red line for
$x_0=-0.1$, and the black line for $x_0=+0.1$.  Note that the value of
$x_1$ is not very important because $Z(s)$ decays very rapidly for large
$\Re(s)$.  Thus, we set $x_1=10$ throughout.  The value of $y_0$ is fixed
at $-0.1$, to avoid integrating over any zeros.  The thick horizontal
line indicates the dimension of the corresponding limit set.  In this
plot, $\theta=20^\circ$.}

\label{fig:run1-20}
\end{figure}

\begin{figure}
\begin{center}
\PSbox{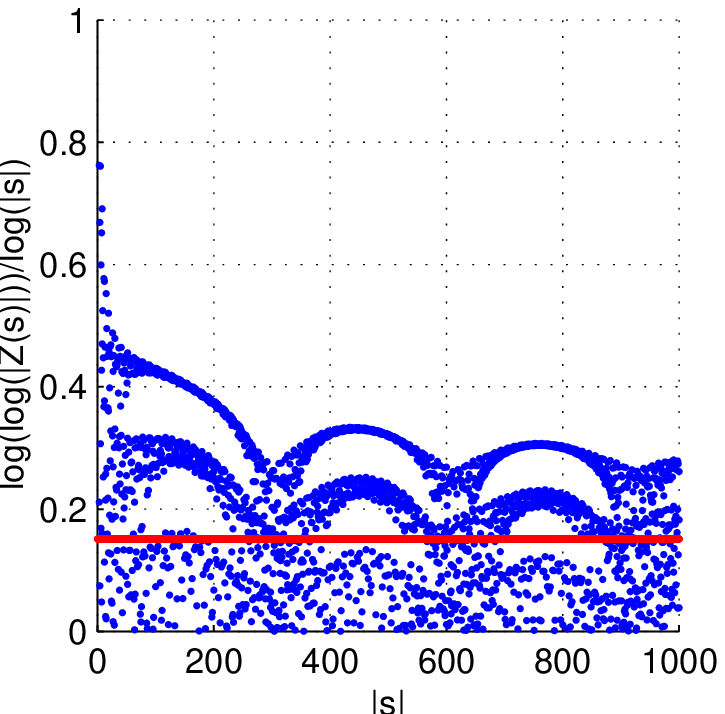}
\end{center}

\caption{This plot shows $\log(\log(|Z(s)|))/\log(|s|)$ for a large
number of points in the rectangle $[-0.2,1]\times[0,10^3]$.  The
horizontal line indicates dimension.  Here, $\theta=20^\circ$.}

\label{fig:zrad-20}
\end{figure}

\begin{figure}
\begin{center}
\PSbox{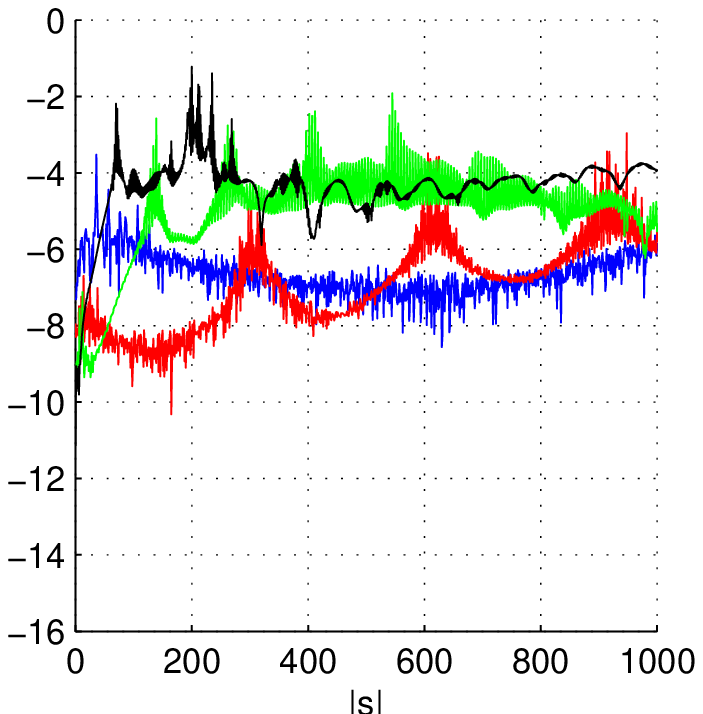}
\end{center}

\caption{Logarithmic plot (base 10) of the modified relative error
$\frac{|R_{12}(s)-R_{13}(s)|}{1+|R_{12}(s)|+|R_{13}(s)|}$ along the line
$\Re(s)=-0.2$, , where $R_N(s)=Z_N'(s)/Z_N(s)$.  The blue curve is $\theta=10^\circ$, the red curve
$\theta=20^\circ$, the green curve $\theta=30^\circ$, and the black
curve $\theta=40^\circ$.}

\label{fig:err-02}
\end{figure}

\begin{figure}
\begin{center}
\PSbox{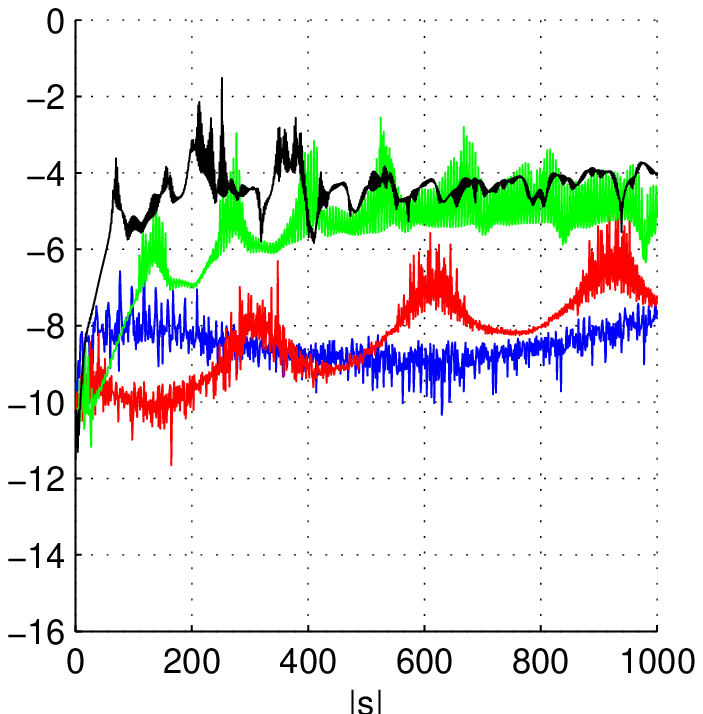}
\end{center}

\caption{Logarithmic plot (base 10) of the modified relative error
$\frac{|R_{12}(s)-R_{13}(s)|}{1+|R_{12}(s)|+|R_{13}(s)|}$ along the line
$\Re(s)=-0.1$,, where $R_N(s)=Z_N'(s)/Z_N(s)$.  The blue curve is $\theta=10^\circ$, the red curve
$\theta=20^\circ$, the green curve $\theta=30^\circ$, and the black
curve $\theta=40^\circ$.}

\label{fig:err-01}
\end{figure}

\begin{figure}
\begin{center}
\PSbox{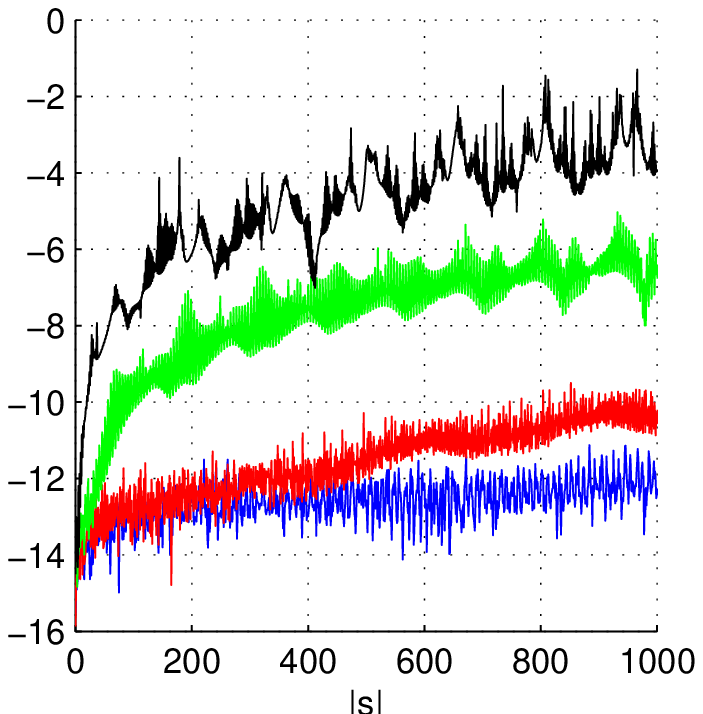}
\end{center}

\caption{Logarithmic plot (base 10) of the modified relative error
$\frac{|R_{12}(s)-R_{13}(s)|}{1+|R_{12}(s)|+|R_{13}(s)|}$ along the line
$\Re(s)=+0.1$, where $R_N(s)=Z_N'(s)/Z_N(s)$.  The blue curve is
$\theta=10^\circ$, the red curve $\theta=20^\circ$, the green curve
$\theta=30^\circ$, and the black curve $\theta=40^\circ$.}

\label{fig:err+01}
\end{figure}

\subsection{Implementation notes.}

To count the number of zeros of $Z(s)$ in a given region $\Omega$ in the
complex plane, we rely on the Argument Principle:
\begin{equation}
  \label{eq:arg-princ}
  |\{s\in\Omega:Z(s)=0\}| = \frac{1}{2\pi{i}}
  \int_{\partial\Omega}{\frac{Z'(s)}{Z(s)} ds} .
\end{equation}
To evaluate $Z(s)$, our main technical tool comes from Jenkinson and
Pollicott {\cite{Jen}}, though we note that the essential ideas were
used in Eckhardt, et. al. {\cite{Eck}} and date back to Ruelle
{\cite{Rue}}.

First, some notation: Let us denote symbolic sequences of length
$|\gamma|=n$ by $\gamma$.  That is,
$\gamma=(\gamma(0),\gamma(1),...,\gamma(n))$, $\gamma(k)\in\{0,1,2\}$,
and $\gamma(0)=\gamma(n)$.  Such sequences represent periodic orbits of
length $n$.  To each such sequence $\gamma$ we associate a composition
of reflections
$\phi_\gamma=\phi_{\gamma(n)}\circ...\circ\phi_{\gamma(1)}:D_{\gamma(0)}
\rightarrow{D_{\gamma(0)}}$.
As $\phi_\gamma$ is a contraction of $D_{\gamma(0)}$ into itself, it has
a unique fixed point $z_\gamma$.

It is shown in Jenkinson and Pollicott that
$Z(s)=\lim_{M\rightarrow\infty}{Z_M(s)}$, where
\begin{equation}
  \label{eqn:jp}
  Z_M(s) = 1 + \sum_{N=1}^{M}{\sum_{r=1}^{N}}{
         \frac{(-1)^r}{r!}\sum_{n\in{P(N,r)}}
         {\prod_{k=1}^{r}{\frac{1}{n_k}{\sum_{|\gamma|=n_k}
         {\frac{|\phi'_\gamma(z_\gamma)|^s}{|1-\phi'_\gamma(z_\gamma)|^2}}}}}},
\end{equation}
where $P(N,r)$ is the set of all $r$-tuples of positive integers
$(n_1,...,n_r)$ such that $n_1+...+n_r=N$.  The series (in $N$)
converges absolutely in $\{s:\Re(s) > -a\}$ for some positive $a$.

Equation (\ref{eqn:jp}) lets us evaluate $Z(s)$ for reasonable values of
$s$ in a straightforward manner.  However, we found two simple but
useful observations during the course of this calculation:
\begin{enumerate}

\item Define
\begin{equation}
   a_n(s) =
  \frac{1}{n}\sum_{|\gamma|=n}{\frac{|\phi'_\gamma(z_\gamma)|^s}
  {|1-\phi'_\gamma(z_\gamma)|^2}}
\end{equation}
and
\begin{equation}
  B_{N,r}(s) = \frac{(-1)^r}{r!}\sum_{n\in{P(N,r)}}
  {\prod_{k=1}^{r}{a_{n_k}(s))}}.
\end{equation}
Then the recursion relation
\begin{equation}
  B_{N,r}(s)=-\frac{1}{r}\sum_{n}{B_{N-n,r-1}(s)\cdot{a_n(s)}},
\end{equation}
with initial conditions $B_{N,1}(s)=a_N(s)$, provides an efficient way
to evaluate the sum in (\ref{eqn:jp}).  A similar relation can be derived
for $Z'(s)$ by differentiation.

\item Recall that the maps $\phi_\gamma$ are compositions of linear
fractional transformations.  Identifying such maps with matrices in
$GL(2,\RR)$ in the usual way, we can compute the numbers
$\phi'_\gamma(z_\gamma)$ via matrix multiplications.  However, such
matrix multiplications can become numerically unstable for larger values
of $|\gamma|$.

An alternative involves the observation that the matrices
$A_\gamma=A_{\gamma(n)}\cdot...\cdot{A_{\gamma(1)}}$ corresponding to
the maps $\phi_\gamma$ have distinct nonzero real eigenvalues.  Let us
denote these eigenvalues by $\lambda_+$ and $\lambda_-$ so that
$|\lambda_+|>|\lambda_-|$.  Then a simple calculation shows that
$\phi'_\gamma(z_\gamma)=\lambda_-/\lambda_+$.  This becomes simply
$(-1)^{|\gamma|}/\lambda_+^2$ if we normalize the determinants of the
generators $A_0$, $A_1$, and $A_2$.  The larger eigenvalue $\lambda_+$
can be easily computed using a naive power method:
\begin{enumerate}

\item
Choose a random $v_0$.

\item
For each $k\geq{0}$, set $v_{k+1}=A_\gamma{v_k}/||A_\gamma{v_k}||$ and
$\lambda_{+}^{(k)}=(A_\gamma{v_k},v_k)$.

\item
Iterate until the sequence $(\lambda_{+}^{(k)})$ converges, up to some
prespecified error tolerance.

\end{enumerate}
The resulting algorithm is slightly less efficient than direct matrix
multiplication, but it is much less susceptible to the effects of
round-off error.

Note that it is certainly possible, even desirable, to apply to this
problem modern linear algebraic techniques, such as those implemented in
ARPACK {\cite{arpack}}.  But, we found that the power method suffices in
these calculations.

\end{enumerate}
These two simple observations lets us calculate the values of $Z(s)$ for
a wide range of values in an efficient manner.  When combined with
adaptive gaussian quadrature, Equation ({\ref{eqn:jp}}) allows us to
evaluate the relevant contour integrals.

\paragraph{Remarks.}
\begin{enumerate}

  \item To calculate the Selberg zeta function $Z_2(s)$ for closed
  geodesics on the quotient space $\Gamma\backslash{\Hh^2}$, we simply sum
  over periodic orbits of even length, and additionally use $a_2(n,s) =
  2a(n,s)$ instead of $a(n,s)$ in the recursion relations above.  This
  counts the number of equivalence classes of orbits correctly.

  \item The work of Pollicott and Rocha {\cite{Pol}} revolves around a
  closely-related trace formula:
  \begin{equation}
    \label{eqn:pr}
    Z(s) = 1 + \sum_{N=1}^{\infty}{\sum_{r=1}^{N}}{
           (-1)^r\sum_{\{[\gamma_1],...,[\gamma_r]\}\in{P_\Gamma(N,r)}}
           {\prod_{k=1}^{r} {\frac{|\phi'_{\gamma_k}(z_{\gamma_k})|^s}
           {|1-\phi'_{\gamma_k}(z_{\gamma_k})|^2}}}}
  \end{equation}
  where $P_\Gamma(N,r)=\{\{[\gamma_1],...,[\gamma_r]\}:
  |\gamma_1|+...+|\gamma_r|=N, \gamma_k\mbox{ primitive}\}$, and
  $[\gamma]$ is the equivalence class of $\gamma$ under shifts.  The
  primary difference between (\ref{eqn:jp}) and (\ref{eqn:pr}) is that
  the latter sums over sets of equivalence classes of primitive periodic
  orbits (equivalent up to shifts), whereas the former sums over all
  periodic orbits.  While it is possible to enumerate primitive periodic
  orbits efficiently,\footnote{For example, by adapting the Sieve of
  Eratosthenes.} Equation (\ref{eqn:jp}) still provides a better
  numerical algorithm, as it is easier to implement and results in
  faster and more stable code.

\end{enumerate}

\subsection{Error analysis.}

Figures {\ref{fig:err-02}}-{\ref{fig:err+01}} show the logarithms (base
10) of the modified relative errors
\begin{equation}
  \frac{|R_{12}(s)-R_{13}(s)|} {1+|R_{12}(s)|+|R_{13}(s)|},
\end{equation}
on the lines $x_0+i[0,10^3]$, for $x_0\in\{-0.2,-0.1,0.1\}$ and where
$R_N(s)=Z'_N(s)/Z_N(s)$.  This formula interpolates between the absolute
and the relative errors, and measures the convergence of the integrand
in (\eqref{eq:arg-princ}).  These results lend some weight to the
reliability (i.e. convergence) of the values of $Z_N'(s)/Z_N(s)$ used in
the calculations above.

\medskip

\noindent
{\sc Acknowledgements.}  The first author would like to
thank J. Anderson for helpful mail discussion. 
In addition to Curt McMullen, the third author would like to 
thank John Lott for introducing him to dynamical zeta functions,
and to Mike Christ for helpful discussions.
KL was supported by a Fannie and John Hertz Foundation
   Fellowship, and  partially  by the Office
   of Science, Office of Advanced Scientific Computing Research,
   Mathematical, Information, and Computational Sciences
   Division, Applied Mathematical Sciences Subprogram, of the
   U.S.\ Department of Energy, under Contract No.\
   DE-AC03-76SF00098.
MZ  gratefully acknowledges the partial support by 
the National Science Foundation under the grant 
DMS-0200732.

\parindent0pt

\end{document}